\documentclass{article}
\usepackage{amsmath, amssymb}
\begin{document}

\author{V. M. Petechuk\footnote{Email: \texttt{vasil-petechuk@rambler.ru}}}
\title{Stability of rings \footnote{The Ukrainian version of this article was published in \cite{Uzhg}.}}
\date{}
\maketitle

\begin{abstract}
The conditions for stability of the elements of linear groups over the associative rings with identity and their connection with the stability of rings are analyzed in the article. The stability of rings which are commutative, satisfy the conditions of stability of rank  $\geq 2$, von Neumann regular, integer-algebraic, nearly local rings introduced by the author, is examined. The most important classical results of H. Bass, L. Vaserstein, S.H. Khlebutin, A.A. Suslin, J.S. Wilson, I.Z.Golubchik, are considered from the unified standpoint.
\end{abstract}

\newtheorem{lr}{Lemma}
\newtheorem{sr}{Corollary}
\newtheorem{vr}{Definition}
\newtheorem{zr}{Remark}
\newtheorem{theo}{Theorem}

\par The study of stability of rings takes root from the well-known commutator formula $SL\left(n,R\right)=\left[SL\left(n,R\right),SL\left(n,R\right)\right]$ and Jordan-Dixon's theorem (1870) about the simplicity of $PSL\left(n,R\right)$ over the field {\it R} when $n\ge 3$. Later, due to the works of Dieudonn\'{e} in the 50s, it turned out that these statements are valid over skew fields as well.

\par The problem of the stability of local rings was formulated and solved in the 60s by Klingenberg. As it turned out, the generalization of the commutator formula over local rings with $n\ge 3$  are the formulas $\left[C\!\left(n,I\right),E\!\left(n,R\right)\right]$ $=E\left(n,I\right)\triangleleft GL\left(n,R\right)$ which hold for all ideals {\it I} of the ring {\it R}; and the generalization of the Jordan-Dixon formula are the inclusions $E\left(n,I_{0} \right)\subset G\subset C\left(n,I_{0} \right)$ for invariant with respect to $E\left(n,R\right)$ subgroups {\it G} of the group $GL\left(n,R\right)$ and their respective ideals $I_{0} $ of the ring {\it R}.

\par Associative rings with identity that satisfy the aforementioned generalizations shall be called {\it stable}.
\par The work of H. Bass in 1964 was the fundamental breakthrough in the proof of stability of wider classes of rings. Therein, he proved the stability of rings that satisfy the condition of stability of rank  $\ge 2$. For instance, semilocal rings that satisfy the rank 1 stability condition are stable as well.

\par In 1972-1977 Wilson, I.Z. Golubchik, A.A. Suslin proved that commutative rings with identity are stable. The stability of rings that are finitely generated over their centers is the generalization of this result. This theorem was obtained independently in the works of A.A. Suslin and L. Vaserstein. Moreover, in 1981 L. Vaserstein has proved the stability of associative rings with identity if their localizations over all maximal ideals of the ring's center satisfy the condition of stability of rank $\ge 2$. In 2001, V.M. Petechuk proved the stability of associative rings with identity should their localizations over all maximal ideals of the ring's center be stable.


\par In 1986 S.H. Khlebutin and L. Vaserstein, independently proved the stability of rings regular in von Neumann sense. In 2001, V.M. Petechuk generalized this result to rings that are integer algebraic over arbitrary subrings of own centers.

\par In 1995, the principally new approach to establishing the stability of weakly Noetherian rings that contain infinite fields in own centers was proposed by I.Z. Golubchik. In 1997 he announced about the stability of the block integer-algebraic rings.

\par Stability of rings is considered in the Chevalley groups as well. Stability of commutative rings in the Chevalley groups is proved in \cite{Pet22, Pet21, Pet23, Pet24} and Lie-type rings in \cite{Pet25}.

\par In the current paper, the stability of associative rings with identity is proved from a wider than earlier notion of the stability of elements of the general linear group.

\par The conditions of stability of elements of linear groups over the associative rings with identity are considered in the paper as well. The stability of rings that satisfy the condition of stability of rank $\ge 2$, commutative, regular in von Neumann sense, block integer-algebraic and the nearly local rings, introduced by the author, follow from these conditions.

\par The statements about the stability of some rings are most systematically laid out in \cite{Pet2, Pet15}. Remarks in the third paragraph on p.122 of \cite{Pet15} do not correspond to reality.

\par Stability of rings plays an important role in applications \cite{Pet26, Pet11, Pet13, Pet3, Pet6, Pet5, Pet7}, particularly in the description of the homomorphisms of linear groups over associative rings \cite{Pet20, Pet8, Pet29, Pet30, Pet15, Pet14}

\par Let {\it R} be an associative ring with identity, $R^{*} $ - group of invertible elements of the ring {\it R}, $J\left(R\right)$ - Jacobson's radical of the ring {\it R},  $\xi R$ and $\xi R^{*} $ - centers of {\it R} and $R^{*}$ respectively.

\par Let $R_{n} = M_{n} \left(R\right)$ denote the ring of all $n\times n$ matrices over {\it R} and $GL\left(n,R\right)=R_{n}^{*} $ - respectively the general linear group of invertible matrices.

\par Under $e_{ij} \in M_{n} \left(R\right)$ we shall understand the matrix with identity at the place $\left(i,j\right)$, and zeros at the rest. Such matrices we shall call {\it standard identity matrices}. The element $t_{ij} \left(r\right)=1+re_{ij} $, where 1 is the identity matrix, $i\ne j,$ 
\\$r\in R$ shall be called a {\it transvection}. Sometimes the identity matrix shall be denoted by {\it E}. If {\it X} is a subset of the ring {\it R} then under $t_{ij} \left(X\right)$ we shall understand the set $\left\{t_{ij} \left(r\right)\left|{\rm \; }r\in X\right. \right\}$ for fixed $i\ne j$ ,  $1\le i,j\le n$, $E_{X} =\left\langle t_{ij} \left(X\right)\right\rangle $ - subgroup of the group $GL\left(n,R\right)$, generated by the set of all $t_{ij} \left(X\right)$,\\
$1\le i,j\le n$.

\par In the particular case when $X=R$ we shall also use a notation $E\!\left(n,R\right)$ $=E_{R}$. Non-identity transvections from $t_{ij} \left(R\right)$ and $t_{ji} \left(R\right)$ shall be called {\it opposite}.

\par Let {\it I} be an arbitrary two-sided ideal of the ring {\it R}. Define by \\
$\Lambda _{I} :R\to {\raise0.7ex\hbox{$ R $}\!\mathord{\left/{\vphantom{R I}}\right.\kern-\nulldelimiterspace}\!\lower0.7ex\hbox{$ I $}} $, $\Lambda _{I} :M_{n} \left(R\right)\to M_{n} \left({\raise0.7ex\hbox{$ R $}\!\mathord{\left/{\vphantom{R I}}\right.\kern-\nulldelimiterspace}\!\lower0.7ex\hbox{$ I $}} \right)$, $\Lambda _{I} :GL\left(n,R\right)\to GL\left(n,{\raise0.7ex\hbox{$ R $}\!\mathord{\left/{\vphantom{R I}}\right.\kern-\nulldelimiterspace}\!\lower0.7ex\hbox{$ I $}} \right)$ the natural homomorphisms of the rings {\it R},  $M_{n} \left(R\right)$ and the group $GL\left(n,R\right)$.

\par Let us define a subgroup $C_{I} =\ker \Lambda _{I} $ that we shall call {\it the main congruence-subgroup of level I} in the group $GL\left(n,R\right)$. The full preimage of the center of the group $GL\left(n,{\raise0.7ex\hbox{$ R $}\!\mathord{\left/{\vphantom{R I}}\right.\kern-\nulldelimiterspace}\!\lower0.7ex\hbox{$ I $}} \right)$ shall be denoted by $C\left(n,I\right)=\Lambda _{I}^{-1} \xi GL\left(n,{\raise0.7ex\hbox{$ R $}\!\mathord{\left/{\vphantom{R I}}\right.\kern-\nulldelimiterspace}\!\lower0.7ex\hbox{$ I $}} \right)$ and under $E\left(n,I\right)$ we shall understand the normal closure of the group $E_{I} $ in $E\left(n,R\right)$. It is easy to see that 
\begin{center}
$\displaystyle E\left(n,I\right)\subseteq C_{I} \subseteq
C\left(n,I\right)$ .
\end{center}

\par Let {\it N} and {\it G} be subgroups of the group $GL\left(n,R\right)$ that are invariant with respect to the group $E\left(n,R\right)$ and {\it N} does not contain non-identity transvections. Under $I_{0}$ we shall understand the largest ideal of ring {\it R} such that $E\left(n,I_{0} \right)\subseteq G$.

\par For elements of an arbitrary group we shall use the notations $a^{b} =bab^{-1}$,  $\left[a,b\right]=aba^{-1} b^{-1} $  and  $\left[a_{1} ,\ldots ,a_{l} \right]=\left[\left[a_{1} ,\ldots ,a_{l-1} \right],a_{l} \right]$ and the commutator formulas
\begin{center}
$\left[ab,c\right]=\left[b,c\right]^{a} \cdot \left[a,c\right]$,           $\left[a,bc\right]=\left[a,b\right]\cdot \left[a,c\right]^{b} $,
\end{center}
and P. Hall's identity
\begin{center}
$\displaystyle \left[a^{-1} ,b,c\right]^{a} \cdot \left[c^{-1} ,a,b\right]^{c} \cdot \left[b^{-1} ,c,a\right]^{b} =1$ .
\end{center}

\par Further  we shall assume that $n\ge 3$.
\par For two non-opposite transvections $t_{ik} \left(x\right)$  and  $t_{lj} \left(y\right)$, $\left(l,j\right)\ne \left(k,i\right)$ the following matrix commutator formulas hold:

\begin{center}$\displaystyle \left[t_{ik} \left(x\right),t_{lj}
\left(y\right)\right]=\left\{\begin{array}{l} {t_{ij} \left(\delta
_{kl} xy\right), i\ne j} \\ {t_{lk} \left(-\delta _{ij} yx\right),l\ne k}
\end{array}\right.,$\end{center}
where $\delta _{ij},\delta _{kl}$ are Kronecker's deltas.

\par From this formula in particular, it follows that the commutator of two non-opposite transvections commutes with each one of them. The following result, first obtained by L. Vaserstein, a shorter proof of which we shall present, holds.
\begin{lr}\label{Petlem1} Let {\it I} be an ideal of the ring {\it R} with identity. Then
$$E\left(n,I\right)=\left\langle t_{ij} \left(I\right)^{t_{ji} \left(R\right)} \left|{\rm \; }1\le i\ne j\le n\right. \right\rangle .$$
\end{lr}

{\it \textbf{Proof}.} Let $T=\left\langle t_{ij} \left(I\right)^{t_{ji} \left(R\right)} \left|{\rm \; }1\le i\ne j\le n\right. \right\rangle $. It is clear that 
\begin{center} $E_{I} \subseteq T\subseteq E\left(n,I\right)$  and  $t_{ij}
\left(I\right)^{t_{kl} \left(R\right)} \subseteq E_{I} $ when
$\left(k,l\right)\ne \left(j,i\right).$\end{center}

\par Therefore, for any transvection $\tau $ the following inclusion is valid
\begin{center} $E_{I}^{\tau } \subseteq T$   and, as a consequence,  $E_{I}^{\tau
\tau _{1} } \subseteq T,$\end{center}
where $\tau _{1} $  is an arbitrary transvection commuting with $\tau$.

\par If the transvection $\tau \notin t_{ij} \left(R\right)$ then the commutators $\left[t_{ji} \left(R\right),\tau \right]$ are transvections, commuting with transvections $\tau $  and  $t_{ji} \left(R\right)$, and it then follows from the matrix commutator formulas that 
\[\displaystyle \left(t_{ij} \left(I\right)^{t_{ji} \left(R\right)} \right)^{\tau } \subseteq \left(t_{ij} \left(I\right)^{\tau } \right)^{t_{ji} \left(R\right)\left[t_{ji} \left(R\right),\tau \right]} \subseteq E_{I}^{t_{ji} \left(R\right)\left[t_{ji} \left(R\right),\tau \right]} \subseteq T.\]

\par When $\tau \in t_{ij} \left(R\right)$, by choosing $s\ne i,j$ we have
\[\begin{aligned}
\left(t_{ij} \left(I\right)^{t_{ji} \left(R\right)}
\right)^{\tau } &\subseteq \left(\left[t_{is} \left(I\right),t_{sj}
\left(R\right)\right]^{t_{ji\left(R\right)} } \right)^{\tau }
\subseteq 
\left[t_{js} \left(I\right)t_{is} \left(I\right),t_{si}
\left(R\right)t_{sj} \left(R\right)\right]^{\tau } \subseteq \\ &\subseteq
\left[E_{I} ,t_{si} \left(R\right)t_{sj}
\left(R\right)\right]\subseteq T.
\end{aligned}\]

\par Therefore, $T^{\tau } \subseteq T$ ,  $E\left(n,I\right)\subseteq T$   and   $E\left(n,I\right)=T$ and so Lemma \ref{Petlem1} is proved.
\begin{sr} Let {\it I}, {\it J}  ideals of the ring {\it R} with identity. Then \\$E\left(n,IJ\right)\subseteq \left[E_{I} ,E_{J}\right].$ In particular,  $E\left(n,I^{2} \right)\subseteq E_{I}.$
\end{sr}

\par Indeed, for any two pairwise distinct numbers $1\le i,j,s\le n$ the following inclusions hold $$t_{ij} \left(IJ\right)^{t_{ji} \left(R\right)} \subseteq \left[t_{is} \left(I\right),t_{sj} \left(J\right)\right]^{t_{ji} \left(R\right)} \subseteq \left[t_{js} \left(I\right)t_{is} \left(I\right),t_{si} \left(J\right)t_{sj} \left(J\right)\right]\subseteq \left[E_{I} ,E_{J} \right].$$

\par In view of Lemma \ref{Petlem1}, $E\left(n,IJ\right)=\left\langle t_{ij} \left(IJ\right)^{t_{ji} \left(R\right)} \left|1\le i\ne j\le n,\right. \right\rangle \subseteq \left[E_{I} ,E_{J} \right].$
By taking $I=J,$ we have $E\left(n,I^{2} \right)\subseteq E_{I}. $

\par The description of the normal structure of linear groups over some rings, as usual, consists of two components. From one side, it is proved that for the group {\it G}, which is normalized by the group $E\left(n,R\right)$, there exists an ideal {\it I} of the ring {\it R} such that 
\begin{center}$\displaystyle E\left(n,I\right)\subseteq
G\subseteq C\left(n,I\right),$
\end{center}
and from the other, the validity of identity $$\left[C\left(n,I\right),E\left(n,R\right)\right]=E\left(n,I\right)\triangleleft GL\left(n,R\right)$$ for an arbitrary ideal {\it I} of the ring {\it R} is proved. It should be noted that both components of the normal structure of linear groups are not always valid at the same time. Moreover, there exist some rings over which none of them is valid. However, the class of rings for which both components of the description of the normal structure of linear groups are valid is quite wide. In particular, it contains commutative rings, rings finitely generated over their centers, and others. The search for conditions for rings, that would be both necessary and sufficient for the aforementioned components of the normal structure of linear groups, is continuing on. The present article represents one such attempt.

\begin{vr} The associative ring {\it R} with identity is called {\it commutator} if for all ideals {\it I} of the ring {\it R} the following identity holds 
 $\left[C\,\left(n,I\right), E\left(n,R\right)\right]=E\left(n,I\right)$  and  $E\left(n,I\right)$ - normal subgroup of  $GL\left(n,R\right)$.
\end{vr}
\begin{vr} The associative ring {\it R} with identity is called weakly-commutator if  there exists a positive integer {\it k} such that  
\begin{center}$\left[C\left(n,I\right),\underbrace{E\left(n,R\right),\ldots ,E\left(n,R\right)}_{k \; times} \right]=E\left(n,I\right)\triangleleft GL\left(n,R\right)$\end{center}  
simultaneously for all ideals {\it I} of the ring {\it R}. The number {\it k} is called the length of the weakly-commutator ring {\it R}.
\end{vr}
\par It is not hard to notice that in an arbitrary associative ring {\it R} with identity one has the inclusion $E\left(n,I\right)\subseteq \left[C\left(n,I\right),E\left(n,R\right)\right]$, where {\it I} is an arbitrary ideal of {\it R}. Moreover, the associative ring {\it R} with identity is commutator if{f} $$\left[C\left(n,I\right),E\left(n,J\right)\right]\subseteq E\left(n,I\right)\bigcap E\left(n,J\right)$$ for all ideals {\it I} and {\it J} of the ring {\it R}. Obviously, commutator rings are weakly-commutator, and in the commutator rings the subgroup $E\left(n,R\right)$ is a normal subgroup of the group $GL\left(n,R\right)$.
\begin{vr} The associative ring {\it R} with identity is called normal if for an arbitrary subgroup {\it G}, invariant with respect to the group $E\left(n, R\right),$ there exists an ideal {\it I} of the ring {\it R} such that
$\displaystyle E\left(n,I\right)\subseteq G\subseteq C\left(n,I\right).$
\end{vr}
\begin{vr} The associative ring {\it R} with identity is called partially normal if an arbitrary subgroup {\it N}, invariant with respect to the group $E\left(n,R\right)$ and not containing non-identity transvections, is contained in $\xi GL\left(n,R\right).$
\end{vr}
\par Obviously the quotient rings of normal rings are partially normal.
\begin{vr} Associative rings that are commutator and normal at the same time are called stable.
\end{vr}
\par It should be highlighted that commutator property, normality and, as a consequence, stability of rings are defined in the group  $GL\left(n,R\right)$ and, therefore, depend on {\it n}.

\begin{lr}\label{Petlem2} The weakly-commutator ring {\it R}, the quotient rings of which are partially normal, is stable.
\end{lr}

{\it \textbf{Proof.}} Let  $I_{0} $ be the largest ideal of {\it R} such that $E\left(n,I_{0} \right)\subseteq G.$ If $\Lambda _{I_{0} } \left(G\right)$ contains transvections, then there exists a nonzero set
\[J_{0} =\left\{ r\in R \left| \Lambda _{I_{0} }
\left(t_{ij} \left(r\right)\right)\in \Lambda _{I_{0} }
\left(G\right) \mbox{for some} \thinspace i \neq j\right.
\right\}.\]

\par It is easy to see that $J_{0}$ is an ideal, containing $I_{0} $. Since $\Lambda _{I_{0} } E\left(n,J_{0} \right)\subseteq \Lambda _{I_{0} } \left(G\right)$, we have $E\left(n,J_{0} \right)\subseteq GC_{I_{0} } $. Therefore, for $r\in J_{0} $ there exists $g\in G$ such that $t_{ij} \left(r\right)g\in C_{I_{0} } \subseteq C\left(n,I_{0} \right)$.
\par Since {\it R} is a weakly-commutator ring of length {\it k}, then 
\begin{center}$\displaystyle \left[C\left(n,I_{0}
\right),\underbrace{E\left(n,R\right),\ldots ,E\left(n,R\right)}_{\mbox{k}
\,\,\mbox{times}} \right]=E\left(n,I_{0} \right)\subseteq
G.$\end{center}
\par Therefore, 
\begin{center} $\left[t_{ij} \left(r\right)g,E\left(n,R\right),\ldots
,E\left(n,R\right)\right]\subseteq G$ and, as a consequence, $t_{ij}
\left(r\right)\subseteq G.$\end{center}

\par This means that $E\left(n,J_{0} \right)\subseteq G$, which contradicts the definition of ideal $I_{0}.$ Therefore, $\Lambda _{I_{0} } \left(G\right)$ does not contain non-identity transvections. Since the ring ${\raise0.7ex\hbox{$ R $}\!\mathord{\left/{\vphantom{R I_{0} }}\right.\kern-\nulldelimiterspace}\!\lower0.7ex\hbox{$ I_{0}  $}} $ is partially normal, we have $\Lambda _{I_{0} } \left(G\right)\subset \xi GL\left(n,{\raise0.7ex\hbox{$ R $}\!\mathord{\left/{\vphantom{R I_{0} }}\right.\kern-\nulldelimiterspace}\!\lower0.7ex\hbox{$ I_{0}  $}} \right)$, i.e. $G\subseteq C\left(n,I_{0} \right).$

\par As a result, 
\begin{center}$\displaystyle E\left(n,I_{0} \right)\subseteq
G\subseteq C\left(n,I_{0} \right).$\end{center}
This proves that {\it R} is a normal ring.
\par Let's prove that {\it R} is a commutator ring.
\par Let $g\in GL\left(n,R\right)$ , $H=E\left(n,R\right)^{g} $  and  $H_{0} =H^{E\left(n,R\right)} $ be a normal closure of {\it H} with respect to $E\left(n,R\right)$. Naturally, $H\subseteq H_{0} $ and $H_{0} $ is an invariant subgroup with respect to the group $E\left(n,R\right)$. Since {\it R} is a normal ring, then there exists an ideal {\it I} of the ring {\it R} such that \begin{center}$\displaystyle E\left(n,I\right)\subseteq H_{0} \subseteq
C\left(n,I\right).$\end{center}

Taking into account the fact that $C\left(n,I\right)$ is a normal subgroup of the group $GL\left(n,R\right)$ and $E\left(n,R\right)^{g}$ $ \subseteq H_{0} \subseteq C\left(n,I\right)$, we have $I=R$,  $E\left(n,R\right)\subseteq H_{0} $ and $\left[E\left(n,R\right),H\right]\subseteq \left[H_{0} ,H\right]$. It follows from $E\left(n,R\right)=\left[E\left(n,R\right),E\left(n,R\right)\right]$ that $H=\left[H,H\right]\subseteq \left[H_{0} ,H\right]$. This proves that \[H_{0} =H^{E\left(n,R\right)} \subseteq \left[E\left(n,R\right),H\right]H\subseteq \left[H_{0} ,H\right].\]

Since the ring {\it R} is weakly-commutator of length {\it k}, it follows that 
\begin{center}$\displaystyle
\left[GL\left(n,R\right),\underbrace{E\left(n,R\right),\ldots
,E\left(n,R\right)}_{k \; times} \right]\subseteq
E\left(n,R\right).$\end{center}

Therefore \begin{center}$\displaystyle H_{0} \subseteq \left[H_{0}
,\underbrace{H,\ldots ,H}_{k \; times} \right]\subseteq
\left[GL\left(n,R\right),\underbrace{H,\ldots ,H}_{k \; times}
\right]\subseteq H.$\end{center}
This means that $H_{0} =H$,  $E\left(n,R\right)\subseteq H$ and $E\left(n,R\right)^{g} \subseteq E\left(n,R\right)$ for all $g\in GL\left(n,R\right)$. Thus, we have proved that $E\left(n,R\right)$ is a normal subgroup of the group $GL\left(n,R\right)$. By taking into account the fact that {\it R} is a weakly-commutator ring of length {\it k}, i.e. $$E\left(n,I\right)=\left[C\left(n,I\right),\underbrace{E\left(n,R\right),\ldots ,E\left(n,R\right)}_{k \; times} \right]$$ we obtain, as a consequence, that $E\left(n,I\right)$ is a normal subgroup of the group $GL\left(n,R\right)$ for all ideals {\it I} of the ring {\it R}.

If $k\ge 2$, we denote
\begin{center}$\displaystyle C_{1}
=\left[C\left(n,I\right),\underbrace{E\left(n,R\right),\ldots
,E\left(n,R\right)}_{k-2 \; times} \right].$\end{center}
Then $\left[C_{1} ,E\left(n,R\right),E\left(n,R\right)\right]=E\left(n,I\right)$  and  $\left[E\left(n,R\right),C_{1} ,E\left(n,R\right)\right]=E\left(n,I\right).$ From P. Hall's commutator identity, by taking account that $E\left(n,I\right)$ is a normal subgroup of the group $GL\left(n,R\right)$, we receive $\left[E\left(n,R\right),E\left(n,R\right),C_{1} \right]\subseteq E\left(n,I\right)$. This means that 
\begin{center}$\displaystyle \left[C_{1} ,E\left(n,R\right)\right]\subseteq
E\left(n,I\right).$\end{center}
\par Therefore, \begin{center}$\displaystyle
\left[C\left(n,I\right),\underbrace{E\left(n,R\right),\ldots
,E\left(n,R\right)}_{k-1 \quad times}
\right]=E\left(n,I\right).$\end{center}
\par Proceeding analogously we obtain $\left[C\left(n,I\right),E\left(n,R\right)\right]=E\left(n,I\right)$. So we proved that {\it R} is a commutator ring and, as a consequence, {\it R} is a stable ring.

From the proof of Lemma \ref{Petlem2} we receive
\begin{sr} Weakly-commutator normal rings are stable.
\end{sr}

It should be noted that in a commutator ring {\it R} for a subgroup {\it L} of the group  $GL\left(n,R\right),$ that for some ideal $I_{0}$ of the ring {\it R} satisfies the condition $E\left(n,I_{0} \right)\subseteq L\subseteq C\left(n,I_{0} \right)$, one has the following inclusions 
\begin{center}$\displaystyle E\left(n,I_{0} \right)\subseteq
\left[L,E\left(n,R\right)\right]\subseteq \left[C\left(n,I_{0}
\right),E\left(n,R\right)\right]=E\left(n,I_{0}
\right).$\end{center}

Therefore, $\left[L,E\left(n,R\right)\right]=E\left(n,I_{0} \right)\subseteq L$. This means that {\it L} is  $E\left(n,R\right)$ -normal subgroup of the group  $GL\left(n,R\right)$  and ideal  $I_{0} $  is uniquely defined by the subgroup {\it L}.
\par Let {\it N} denote an $E\left(n,R\right)$-invariant subgroup of the group $GL\left(n,R\right),$ which does not contain non-identity transvections.

If {\it I} is a two-sided ideal of the ring {\it R} then the annihilator $$AnnI=\left\{r\in R\left|rI=Ir=0\right. \right\}$$ of ideal {\it I} in {\it R} is a two-sided ideal as well.

\begin{lr}\label{Petlem3} Let {\it R} be an associative ring with identity, $g=\left(g_{ij} \right)\in N$ and there exists $x\in R$ such that $g_{ij} x=0$ for some fixed $1\le i,j\le n$. Then \\$g\in C\left(n,AnnRxR\right)$ if $i\ne j$  and  $x=0$ otherwise.
\end{lr} 

{\it\textbf{Proof.}} It is not hard to see that in the case $k\ne j$ the {\it i}th row of the matrix $g_{1} =\left[g,t_{jk}
\left(x\right)\right]\in N$ coincides with the {\it i}th row of the matrix $t_{jk} \left(-x\right)$.

Suppose that $i\ne j$. Then the {\it i}th row of the matrix $t_{jk} \left(-x\right)$ and the identity matrix coincide. If $g_{1} \ne 1$, then {\it N} contains transvections of the type $\left[g_{1} ,t_{li} \left(R\right)\right]$ for all $l\ne i$ ,  $1\le l\le n$. Since {\it N} does not contain non-identity transvections,
 we have $g_{1} =1$ for all $k\ne j$. It follows from the identity $\left[g,t_{jk} \left(x\right)\right]=1$ that $xg_{ks} =0$ for all $s\ne k,$ $g_{sj} x=0$ for all $s\ne j$ and $xg_{kk} =g_{jj} x.$

Following analogously we prove that $g^{-1} $ commutes with all matrices $t_{mk} \left(x\right)$, where $m\ne k$,  $x\in R$. Consequently, {\it g} commutes with all the transvections from the group $E_{x} $. This is equivalent to the condition that $gx=xg$ is a scalar matrix (however, not necessarily central).

Since $g_{ij} xR=0$, we have that $gxr=xrg$ is a scalar matrix for all $r\in R.$ Taking into account $i\ne j$ we have $xRg_{ij} =0$ and, consequently, $RxRg_{ij} =0.$ As above we prove that $r'xrg=gr'xr$ is a scalar matrix for all {\it r},  $r'\in R.$ This means that elements from $RxR$ annihilate from the left and right the elements $g_{pq},$   $g_{pp} -g_{pq} $ of matrix {\it g} for all $1\le p\ne q\le n.$ Hence, it is proved that $g\in C\left(n,AnnRxR\right).$

In particular, when {\it g} has a zero non-diagonal element which, obviously, is annihilated by all the elements of the ring $R,$ we obtain  $g\in C\left(n, AnnR\right)=\xi GL\left(n,R\right)$.

\par Since $g_{1} =\left[g,t_{jk} \left(x\right)\right]\in N$  and  $g_{1} $ has a zero non-diagonal element(as $n\ge 3$), we have  $\left[g,t_{jk} \left(x\right)\right]\in \xi GL\left(n,R\right).$

Let us consider the case when $i=j$.  As $\left[g,t_{ik} \left(x\right)\right]\in \xi GL\left(n,R\right)$ and the {\it i}th row of matrix $gxe_{ik} g^{-1} $ is all zero, then for some element $r\in \xi R\bigcap R^{*} $, the {\it i}th row of matrix $rt_{ik} \left(x\right)-E=\left(r-1\right)E+rxe_{ik}$ is all zero as well. Therefore, $r=1$  and  $x=0$.
\begin{zr} Following analogously, the arguments in Lemma \ref{Petlem3} remain valid when instead of the equality $g_{ij} x=0$ one considers the equality $xg_{ij} =0$ .
\end{zr}

It follows from Lemma \ref{Petlem3} that the diagonal elements of matrices of the group {\it N} do not have left or right zero dividers. In particular, the diagonal elements of the matrices $g\in N$ cannot be equal to zero.

\begin{sr}\label{Petnas3} If $g\in N$ and for some $x\in R$ the commutator $\left[g,t_{ij} \left(x\right)\right]$ has a zero element, then $g\in C\left(n,AnnRxR\right)$.
\end{sr}

{\it\textbf{Proof.}} Since {\it N} is an $E\left(n,R\right)$-invariant group, then $\left[g,t_{ij} \left(x\right)\right]\in N$ and, according to Lemma 3, the inclusion $\left[g,t_{ij} \left(x\right)\right]\in \xi GL\left(n,R\right)$ holds. Therefore, there exists $r\in \xi R\bigcap R^{*} $ such that $gt_{ij} \left(x\right)g^{-1} =rt_{ij} \left(x\right)$. This means that $gxe_{ij} =\left(rt_{ij} \left(x\right)-E\right)g=\left(r-1\right)g+rxe_{ij} g$. In such a case $\left(r-1\right)g_{ll} =0$ where $l\ne i,j$. Since the diagonal elements of matrix {\it g} do not have zero divisors, then $r=1$ and $gxe_{ij} =xe_{ij} g$. Thus, $g_{si} x=0$  for all $s\ne i$. According to Lemma \ref{Petlem3}, the inclusion $g\in C\left(n,AnnRxR\right)$ holds.

\begin{lr}\label{Petlem4} Let $g\in N$  and  $x_{1} ,\ldots ,x_{n} \in R$ such that $g_{i1} x_{1} +\cdots +g_{in} x_{n} =0$ and at least one of the elements $x_{i} $ is equal to zero. Then $$g\in C\left(n,Ann\left(Rx_{1} R+ \ldots +Rx_{n} R\right)\right).$$
\end{lr}

{\it\textbf{Proof.}} Suppose that $x_{j} =0$ for some $1\le j\le n.$ Then the {\it i}th row of the commutator $g_{1} =\left[g,t_{1j} \left(x_{1} \right)\cdots t_{nj} \left(x_{n} \right)\right]\in N$ coincides with the {\it i}th row of the matrix $t_{ij} \left(-x_{i} \right)$. Since it contains zero non-diagonal elements, then, according
 to Lemma 3, the commutator $g_{1} \in \xi GL\left(n,R\right)$. Recalling that there is an identity on the $\left(i,i\right)$th place in matrix $t_{ij} \left(-x_{i} \right)$, we have $g_1=1$. Therefore, for an arbitrary $1\le l\le n$ there exists $1\le s\ne k\le n$ such that $x_{l} g_{sk} =0$. According to Lemma 3, $g\in C\left(n,AnnRx_{l} R\right)$ for all $1\le l\le n$. In this case, the elements $g_{pq} $ and $g_{pp} -g_{qq} $, where $1\le p\ne q\le n$ are annihilated from left and right by elements of the ideals $Rx_{l} R$ for all $1\le l\le n$ and, therefore, by elements of the sum $Rx_{1} R+\cdots +Rx_{n} R.$ This proves the inclusion \\$g\in C\left(n,Ann\left(Rx_{1} R+\cdots +Rx_{n} R\right)\right).$

\begin{zr}\label{Petzauv2} Analogously, if $g\in N$  and  $x_{1},\ldots ,x_{n} \in R$ such that \\
$x_1g_{1j}+\ldots+x_n g_{nj}=0$ and at least one of the elements $x_{i} $ is zero, then \\$g\in C\left(n,Ann\left(Rx_{1} R+\cdots + \right.\right.$ $\left.\left.Rx_{n} R\right)\right).$
\end{zr}
\begin{sr}\label{Petnas4} If $g\in N$ and at least one, not necessarily diagonal, element of the matrix {\it g} has a left or right inverse, then $g\in \xi GL\left(n,R\right)$.
\end{sr}

{\it\textbf{Proof.}} Suppose that $g_{ij}^{-1} $ is a left inverse of $g_{ij} $. For $k\ne i$ put \\$x_{i} =-g_{kj}g_{ij}^{-1}$,  $x_{k} =1$ and $x_{s} =0$ for all $s\ne i,k$. Then $x_1g_{1j}+\ldots+x_n g_{nj}=0$ and $Rx_{1} R+\cdots +Rx_{n} R=R$. According to Remark \ref{Petzauv2} of Lemma \ref{Petlem4}, the inclusion $g\in \xi GL\left(n,R\right)$ is valid.

\begin{lr}\label{Petlem5} Let $g=g_{1} g_{2} \in N$ be such that $\left(g_{1} \right)_{ki} =\delta _{ki} $ and $\left(g_{2} \right)_{kj} =\delta _{kj} $ for all $1\le k\le n$ and fixed $1\le i,j\le n$. Then $g\in \xi GL\left(n,R\right)$.
\end{lr}

{\it\textbf{Proof.}} The condition of Lemma \ref{Petlem5} means that the {\it i}th column of the matrix $g_{1}-1$ and the {\it j}th column of matrix $g_{2} -1$ are zero. If $i=j$, then the {\it i}th column of matrix {\it g} coincides with the {\it i}th column of the identity matrix. Therefore, {\it g} has an invertible element and, in correspondence with Corollary \ref{Petnas4} of Lemma \ref{Petlem4}, $g\in \xi GL\left(n,R\right)$.
\par Suppose that $i\ne j$. Without loss of generality, up to the similarity by a matrix from the group $E\left(n,R\right)$, we can assume that $i=1$ and $j=n$. Then
\begin{center}$\displaystyle g=\left(\begin{array}{cc} {1} & {x} \\ {0} & {A}
\end{array}\right)\left(\begin{array}{cc} {B} & {0} \\ {y} & {1}
\end{array}\right),$ \end{center}
where $A,B\in GL\left(n-1,R\right)$ and {\it x}{\it , }{\it y} - rows of length {\it n}{\it -}1. Let
\begin{center}
 $X=\left(\begin{array}{cc} {1} & {-xA^{-1} } \\ {0} & {E} \end{array}\right)$,        $Y=\left(\begin{array}{cc} {E} & {0} \\ {-y} & {1} \end{array}\right)$,       $g_{0} =XgY$
\end{center}
where {\it E} is an identity $\left(n-1\right)\times \left(n-1\right)$ matrix. Then 
\begin{center}
$\displaystyle g_{0} =\left(\begin{array}{cc} {1} & {0} \\ {0} & {A} \end{array}\right)\left(\begin{array}{cc} {B} & {0} \\ {0} & {1} \end{array}\right)$.
\end{center}
Since $X\in E\left(n,R\right)$  and {\it Y} commutes with $t_{n1} \left(1\right)$, then $X\left[g,t_{n1} \left(1\right)\right]X^{-1} \in N$ and $\left[g_{0} ,t_{n1} \left(1\right)\right]\in N\left[X,t_{n1} \left(1\right)\right]$. On the other hand, the first row of matrices 
\begin{center} $\left(\begin{array}{cc} {1} & {0} \\ {0} & {A}
\end{array}\right)\left[\left(\begin{array}{cc} {B} & {0} \\ {0} &
{1} \end{array}\right),t_{n1}
\left(1\right)\right]\left(\begin{array}{cc} {1} & {0} \\ {0} & {A}
\end{array}\right)^{-1}$,           $\left[\left(\begin{array}{cc}
{1} & {0} \\ {0} & {A} \end{array}\right),t_{n1}
\left(1\right)\right]$
\end{center}
and, as a consequence, of matrix $\left[g_{0} ,t_{n1} \left(1\right)\right]$ coincides with the first row of the identity matrix. Hence,
\begin{center}
$\displaystyle \left[g_{0} ,t_{n1} \left(1\right),t_{n1} \left(1\right),t_{n1} \left(1\right)\right]=E$.
\end{center}
This means that 
\begin{center} $\displaystyle h=\left[X,t_{n1} \left(1\right),t_{n1}
\left(1\right),t_{n1} \left(1\right)\right]\in N$.
\end{center}
Since the second row of matrix {\it h} coincides with the second row of the identity matrix, then, in accordance with the corollary \ref{Petnas4} of Lemma \ref{Petlem4}, $h=E$.

Let $-xA^{-1} =\left(x_{2} ,\ldots ,x_{n} \right)$. Direct computation shows that from the equality $h=E$ one has $x_{n}^{4} =0$. Thus, $\left[t_{n1} \left(1\right),X\right]_{11} =1-x_{n} \in R^{*} $.
\par As a corollary we obtain that element of the matrix $\left[g_{0} ,t_{n1} \left(1\right)\right]\cdot \left[t_{n1} \left(1\right),X\right],$ which is at the position $\left(1,1\right),$ coincides with $1-x_{n} $ and, therefore, is invertible.

Taking into account the fact that {\it Y} commutes with $t_{n1} \left(1\right),$ from the commutator formulas one has an equality 
\begin{center}
$\displaystyle \left[g,t_{n1} \left(1\right)\right]=\left[X^{-1} g_{0} ,t_{n1} \left(1\right)\right]=\left[g_{0} ,t_{n1} \left(1\right)\right]^{X^{-1} } \left[X^{-1} ,t_{n1} \left(1\right)\right]$.
\end{center}
	Hence, its obvious corollary has place
\begin{center} $\displaystyle \left[g,t_{n1}
\left(1\right)\right]^{X} =\left[g_{0} ,t_{n1}
\left(1\right)\right]\left[t_{n1} \left(1\right),X\right]$.
\end{center}
Since $g\in N$, then $\left[g,t_{n1} \left(1\right)\right]^{X} \in N$ and, according to the Corollary \ref{Petnas4} of Lemma \ref{Petlem4}, $\left[g,t_{n1} \left(1\right)\right]^{X} \in \xi GL\left(n,R\right)$ and, as a consequence, $\left[g,t_{n1} \left(1\right)\right]\in \xi GL\left(n,R\right)$ and $g\in \xi GL\left(n,R\right).$

We shall need some statements valid for arbitrary rings.
The following has place
\begin{lr}\label{Petlem6} Let {\it a}{\it , }{\it b}{\it , }{\it c} be some elements of the associative ring {\it R} with identity. Element $1+ab\in R^{*}$ if and only if $1+ba\in R^{*} $. In particular, if $a^2=b^2=0$ and $1+ab\in R^{*}$ then one has a decomposition
\begin{center}
$\displaystyle 1+ab=\left(1+b\left(1-\gamma \right)\right)\left[1-b,1+a\right]\left(1+\left(1-\gamma \right)a\right)\left(1+ba\right)$,
\end{center}
where  $\gamma =\left(1+ab\right)^{-1}.$
\end{lr}

{\it\textbf{Proof.}} The first half of Lemma \ref{Petlem6} follows from the fact that the equality $(1+ab)c=1$ draws the equality $\left(1+ba\right)\left(1-bca\right)=1$.
\par If $a^2=b^2=0$ and  $1+ab\in R^{*}$, then $1+a$,  $1-b$ - invertible elements and $\gamma \left(1+ab\right)=\left(1+ab\right)\gamma =1,$ where $\gamma\!=\!\left(1+ab\right)^{-1}.$ Thus, $\gamma ab\!+\!\gamma\!-\!1=ab\gamma+\gamma-1=0.$

As a corollary we obtain \\$\displaystyle b\gamma aba+b\gamma a-ba=b\gamma ^{2} a-b\gamma
a+b\gamma ab\gamma a=b\gamma ab+b\gamma -b=ab\gamma a+\gamma a-a=0$
 and also  $a\left(1-\gamma \right)=\left(1-\gamma \right)b=0$
and  $a\gamma a=b\gamma b=0$.
\par In this case
\[\begin{aligned}
\left(1-a\right)\left(1+\left(1-\gamma \right)a\right)&=1-\gamma
a-a\left(1-\gamma \right)a=1-\gamma a \quad \mbox{and}\\
\left(1+b\left(1-\gamma \right)\right)\left(1-b\right)&=1-b\gamma
-b\left(1-\gamma \right)b=1-b\gamma. \end{aligned} \]

\par  By direct calculations we establish that
\[\begin{aligned}
\left(1-b\gamma\right)\left(1+a\right)\left(1+b\right)\left(1-\gamma
a\right)&=\left(1+a-b\gamma -b\gamma a\right)\left(1+b-\gamma a-b\gamma a\right)= \\ &=1+ab-b\gamma a 
\end{aligned}\]
and $1+ab=\left(1+ab-b\gamma a\right)\left(1+ba\right)$.
This means that the following equality holds 
\begin{center}$\displaystyle 1+ab=\left(1-b\gamma
\right)\left(1+a\right)\left(1+b\right)\left(1-\gamma
a\right)\left(1+ba\right),$
\end{center} from which the formula of Lemma \ref{Petlem6} follows.

If in Lemma \ref{Petlem6} one puts $a=xe_{ik}$  and $b=-ye_{lj}$, where $x, y \in R,$ then in the case $i\neq k$ and $l\neq j,$ the equalities $a^2=b^2=0$ and $ab=-\delta _{kl} xye_{ij},$  $ba=-\delta _{ij} yxe_{lk}$ are valid.
\par If $i\ne j$, then $ba=0$ and $\gamma =\left(1+ab\right)^{-1} =1-ab$. Thus, \\$1+b\left(1-\gamma \right)=1+\left(1-\gamma \right)a=1$. Taking into account that $1+a=t_{ik} \left(x\right)$,  $1-b=t_{lj} \left(y\right)$,  $1+ab=t_{ij} \left(-\delta _{kl} xy\right)$ we receive a formula $t_{ij} \left(-\delta _{kl} xy\right)=\left[t_{lj} \left(y\right),t_{ik} \left(x\right)\right]$ from which one obtains well-known commutator formulas
\begin{center} $t_{ij} \left(\delta _{kl} xy\right)=\left[t_{ik} \left(x\right),t_{lj} \left(y\right)\right]$ and $t_{lk}\left(\delta _{ij} yx\right)=\left[t_{lj}(y), t_{ik}(x)\right].$\end{center}

Let $g\in GL\left(n,R\right)$  and  $U=ge_{ii}$, $V=e_{ij} g^{-1}$, $r,r_{0} \in R$. Obviously, $VU=\delta _{ij} e_{ii}$. Further we shall assume that $i\ne j$. Hence, $VU=0$.
\par Let $x_{1} ,\ldots ,x_{n} $ be some elements of {\it R} and $V_{0} =x_{1} e_{i1} +\cdots +x_{n} e_{in}$, \\$\alpha =x_{1} g_{1i} +\cdots +x_{n} g_{ni}$. Suppose that $x_{l} =0$ for some $1\le l\le n$ and $x_{k} =rg_{jk}^{-1} $ for some $1\le k\le n$. Put $W=rV-V_{0}$. Obviously, $V_{0} U=\alpha e_{ii} =-WU$ and $W_{ik} =0$.
\par Under these notations,
\begin{center}
$\displaystyle t_{ij} \left(r_{0} r\right)^{g} =1+ge_{ii} r_{0} re_{ij} g^{-1} =1+Ur_{0} rV$.
\end{center}
It is esay to check that 
\[\begin{aligned}
t_{ij} \left(r_{0} r\right)^{g} \left(1-Ur_{0}W\right)&=\left(1+Ur_{0} rV\right)\left(1-Ur_{0} W\right)= \\ &=1+Ur_{0}
\left(rV-W\right)=1+Ur_{0} V_{0}.
\end{aligned}\]

Suppose that $1+Ur_{0} V_{0} \in GL\left(n,R\right)$. Then $1-Ur_{0} W\in GL\left(n,R\right)$ and 
\begin{center}$t_{ij} \left(r_{0} r\right)^{g} =\left(1+Ur_{0} V_{0} \right)\left(1-Ur_{0} W\right)^{-1}$\end{center} 
and, in particular, $t_{ij} \left(r\right)^{g} =\left(1+UV_{0} \right)\left(1-UW\right)^{-1}.$
\par The representation of the matrix $t_{ij} \left(r\right)^{g} $ thus obtained is useful because the {\it i}th rows of each of the matrix $V_{0} $ and {\it W} contain at least one zero element. In the long run, this will allow to use Lemma \ref{Petlem6} and decompose the respective commutators into the product of transvections and diagonal elements. 
\par For $g' \in GL(n,R)$ we analogously define $U'=g'e_{ii},$ $V'=e_{ij}(g')^{-1},$ \\$V'_0=x'_1e_{i1}+\ldots+x'_ne_{in},$
$\alpha'=x'_1g'_{1i}+\ldots+x'_ng'_{ni},$ where $x'_1,\ldots,x'_n \in R,$ $x'_l=0$ and $x'_k=r(g')^{-1}_{jk}$ for corresponding $1\leq i,j,l,k\leq n$ and $r_0,r \in R.$ Let $W'=rV'-V'_0$ and $1+U'r_0V'_0 \in GL(n,R).$ Then 
\begin{center}$t_{ij} \left(r_{0} r\right)^{g'}
=\left(1+U'r_{0} V'_{0} \right)\left(1-U'r_{0} W'\right)^{-1}.
$\end{center}
One has 
\begin{lr}\label{Petlem7}{\bf(main)} Let {\it I} and {\it J} be ideals of {\it R}, $c\in \xi R$,  $r\in R$,  $1+V_{0}c^{2}  U\in GL\left(n,R\right),$ $1+\left(V_0-x_ke_{ik}\right)c^2U\in GL(n,R),$ $1+V'_{0}c^{2} U'\in GL\left(n,R\right),$\\
$1+\left(V'_0-x'_ke_{ik}\right)c^2U'\in GL(n,R),$ $h\in C_{I}, $ $g'=gh^{-1},$  $x_{s} -x'_{s} \in I,$ $1\le s\le n$ (for  $s=l,k$ it is automatically satisfied due to $x_{l} =x'_{l} =0$, $x_{k}=rg_{jk}^{-1} $, $x'_{k} =r\left(g'_{jk} \right)^{-1} $). Then

1)  $\left[h,t_{ij} \left(c^{2} r\right)\right]^{g'} \subseteq E\left(n,cI\right)\bigcap E\left(n,J\right)$, if  $r\in J$;

2)  $g\in C\left(n,AnnRc^{2} rR\right)$, if  $g\in N.$
\end{lr}

{\it\textbf{Proof.}} Assume that $u_{s} =g_{si}$, $d_{s} =1+\alpha c^{2} e_{ss}$, where $1\le s\le n$. Since $V_{0} U=\alpha e_{ii}$ and by the condition $1+\alpha c^{2} e_{ii} =1+V_{0}c^{2}  U \in GL\left(n,R\right)$, then $1+\alpha c^{2} \in R^{*} $ and
\begin{center}
$\displaystyle d_{s} =diag\left(1,\ldots ,\underbrace{1+\alpha c^{2}
}_{\mbox{s-th place}} ,1,\ldots ,1\right)\in GL\left(n,R\right)$ .
\center \end{center}
In particular, $d_{i} =1+V_{0} c^{2} U.$
\par Analogously one can define $u'_{s}=g'_{si},$ $d'_s=1+\alpha'c^2e_{ss},$ $1+\alpha'c^{2} \in R^{*}$, $d'_{s} \in GL\left(n,R\right)$, $d'_{i} =1+V'_{0} c^{2} U'.$
Let us consider the matrices $a=\left(Ue_{il} -u_{l} e_{ll} \right)c$  and  $b=e_{li} V_{0} c.$ Obviously,
\begin{center}
 $a=\left(\begin{array}{ccc} {0} & {\begin{array}{l} {u_{1} } \\ {\vdots } \\ {u_{l-1} } \end{array}} & {0} \\ {0} & {0} & {0} \\ {0} & {\begin{array}{l} {u_{l+1} } \\ {\vdots } \\ {u_{n} } \end{array}} & {0} \end{array}\right)c$        and         $b=\left(\begin{array}{ccc} {0} & {0} & {0} \\ {x_{1} \cdots x_{l-1} } & {0} & {x_{l+1} \cdots x_{n} } \\ {0} & {0} & {0}
 \end{array}\right)c$\end{center}

It is not hard to see that $aRa=bRb=e_{ll} Ra=bRe_{ll} =0$. From the definition of matrices {\it a} and {\it b} it follows that 
\begin{center}
 $a+u_{l} e_{ll} c=Ue_{il} c$,        $e_{il} b=V_{0} c$,       $e_{ll} b=b$,        $be_{ll} =0$,         $bd_{l} =b$.
\end{center}
Moreover, {\it ab} is a matrix with zero {\it l}-th row and {\it l}-th column, $ba=\alpha c^{2} e_{ll} $,   $1+ba=1+\alpha c^{2} e_{ll} =d_{l} \in GL\left(n,R\right)$.
\par In accordance with Lemma \ref{Petlem6}, $1+ab\in GL\left(n,R\right)$ and
\begin{center} $1+ab=\left(1+b\left(1-\gamma\right)\right)\left[1-b,1+a\right]\left(1+\left(1-\gamma
\right)a\right)\left(1+ba\right)$, \\where  $\gamma=\left(1+ab\right)^{-1}.$ \end{center}
Obviously, $1-\gamma $ is a matrix with zero {\it l}-th row and {\it l}-th column and matrices $1\pm a$ ,  $1\pm b$ ,  $1+d_{l} u_{l} cb$ ,  $1+\left(1-\gamma \right)a$ ,  $1+b\left(1-\gamma \right)$ are products of transvections.

\par Consider an equality 
$$1+Uc^{2} V_{0} =1+Uce_{il} b=1+\left(a+u_{l} e_{ll} c\right)b=\left(1+ab\right)\left(1+u_{l} cb\right)$$
and the equality \begin{center}$\left(1+ba\right)\left(1+u_{l}
cb\right)=d_{l} \left(1+u_{l} cb\right)=\left(1+d_{l} u_{l}
cb\right)d_{l}.$\end{center}
Let us denote 
\begin{center}$T_l(V_0)=\left(1+b\left(1-\gamma
\right)\right)\left[1-b,1+a\right]\left(1+\left(1-\gamma
\right)a\right)\left(1+d_{l} u_{l} cb\right).$\end{center}
Clearly, $T_l\left(V_0\right)$ is a product of transvections contained in the group $E_{cR} $ and

\[ 1+Uc^{2} V_{0}=\left(1+ab\right)\left(1+u_{l} cb\right)= \\\]
\[=\left(1+b\left(1-\gamma\right)\right)\left[1-b,1+a\right]\left(1+\left(1-\gamma\right)a\right)\left(1+ba\right)\left(1+u_{l}
cb\right)=T_{l}\left(V_0\right)d_{l}.  \]
It should be noted that the construction of decomposition of matrices $1+Uc^{2} V_{0}$ into a product of transvections and diagonal elements is valid for an arbitrary column $U=\left(u_{1} \ldots u_{n} \right)^{T} $ and arbitrary row $V_{0} =\left(x_{1} \ldots x_{n} \right)$ for which $x_{l} =0$ and $1+V_{0} Uc^{2} \in R^{*}$.
Its view is determined by the formula in Lemma \ref{Petlem6}.
\par Thus $1-Uc^2W=T_k(-W)d_k \in GL(n,R).$ Since \\$t_{ij}(rc^2)^g=(1+Uc^2V_0)(1-Uc^2W)^{-1},$ then 
\begin{center}$t_{ij}(rc^2)^g=T(g)d_ld_k^{-1},$ where $T(g)=T_{l}(V_0)(T_k(-W)^{-1})^{d_ld_k^{-1}}.$\end{center}
By taking into account analogous arguments, we have 
\begin{center}              $t_{ij} \left(rc^{2} \right)^{g'}
=T\left(g'\right)d'_{l} \left(d'_{k} \right)^{-1}, $\end{center}
where $T\left(g'\right)$ is exactly the same product of transvections modulo the ideal $cI$ as $T\left(g\right)$. Therefore, 
\begin{center}
$\displaystyle \left[h,t_{ij} \left(rc^{2} \right)\right]^{g'}
=t_{ij} \left(rc^{2} \right)^{g} t_{ij} \left(-rc^{2} \right)^{g'}
=T\left(g\right)d_{l} d_{k}^{-1} d'_{k}\left(d'_{l} \right)^{-1}
T\left(g'\right)^{-1}.$ \end{center}
It is not hard to see that in the case $x\in R^{*} $ one has the formula
\begin{center}$\left(\begin{array}{cc} {x} & {0} \\ {0} & {x^{-1} }
\end{array}\right)=\left(\begin{array}{cc} {1} & {x} \\ {0} & {1}
\end{array}\right)\left(\begin{array}{cc} {1} & {0} \\ {-x^{-1} } &
{1} \end{array}\right)\left(\begin{array}{cc} {1} & {x} \\ {0} & {1}
\end{array}\right)\left(\begin{array}{cc} {1} & {-1} \\ {0} & {1}
\end{array}\right)\left(\begin{array}{cc} {1} & {0} \\ {1} & {1}
\end{array}\right)\left(\begin{array}{cc} {1} & {-1} \\ {0} & {1}
\end{array}\right)$\end{center}
In particular, if {\it I} is an ideal of the ring {\it R} and $x\in 1+cI$, then the matrix
\begin{center}$\left(\begin{array}{cc} {x} & {0} \\ {0} & {x^{-1} }
\end{array}\right)^{\left(\begin{array}{cc} {1}
& {-1} \\ {0} & {1} \end{array}\right)}=$\end{center} \begin{center}$=\left(\begin{array}{cc} {1}
& {x\!-\!1} \\ {0} & {1} \end{array}\right)\left(\begin{array}{cc}
{1} & {0} \\ {1\!-\!x^{-1} } & {1}
\end{array}\right)\left(\begin{array}{cc} {1} & {1\!-\!x} \\ {0} & {0}
\end{array}\right)\left[\left(\begin{array}{cc} {1} & {x\!-\!1} \\ {0} &
{1} \end{array}\right)\!,\!\left(\begin{array}{cc} {1} & {0} \\ {-1}
& {1} \end{array}\right)\right]$\end{center}
is contained in $E(2,cI).$

If $1+\alpha'c \in R^{*}$ and  $\alpha -\alpha '\in I$ then $\left(1+\alpha c\right)\left(1+\alpha 'c\right)^{-1} \in 1+cI$. Hence, 
$$\begin{array}{l} {\left(\begin{array}{cc} {1+\alpha c} & {0} \\ {0} & {\left(1+\alpha c\right)^{-1} } \end{array}\right)\left(\begin{array}{cc} {1+\alpha 'c} & {0} \\ {0} & {\left(1+\alpha 'c\right)^{-1} } \end{array}\right)^{-1} =} \\ {=\left(\begin{array}{cc} {\left(1+\alpha c\right)\left(1+\alpha 'c\right)^{-1} } & {0} \\ {0} & {\left(1+\alpha c\right)^{-1} \left(1+\alpha 'c\right)} \end{array}\right)\in E\left(2,cI\right)}. \end{array}$$

This means that $d_{l} d_{k}^{-1} d'_{k} \left(d'_{l} \right)^{-1} \in E\left(n,c^{2} I\right).$
\par Taking into account that for arbitrary transvections $\tau _{1} $  and  $\tau _{2} $ from $E_{cR} $ such that $\tau _{1} \equiv \tau _{2} \mod E_{cI},$ one has the inclusions $$\tau _{1} E\left(n,cI\right)\tau _{2}^{-1} \subseteq E\left(n,cI\right),$$ we obtain $$T\left(g\right)E\left(n,cI\right)T\left(g'\right)^{-1} \subseteq E\left(n,cI\right).$$
Thus, it is proved that 
\begin{center}
$\displaystyle \left[h,t_{ij} \left(rc^{2} \right)\right]^{g'}
\subseteq T\left(g\right)E\left(n,cI\right)T\left(g'\right)^{-1}
\subseteq E\left(n,cI\right).$\end{center}
Let {\it J} be an ideal of {\it R}, $r\in J$. Then $V_{0} +W=rV$ is a matrix over {\it J} and  $x_{k} \in J$.

Let us prove that $t_{ij} \left(r\right)^{g} \in E\left(n,J\right)$.

As above,
\begin{center}
$\displaystyle t_{ij} \left(rc^{2} \right)^{g} =\left(1+Uc^{2} V_{0}
\right)\left(1-Uc^{2} W\right)^{-1} =T\left(g\right)d_{l} d_{k}^{-1}
$ \end{center} is a product of transvections and diagonal elements, which, as is known, normalize the group $E(n,I).$

Consider $V_{0}^{*} =V_{0} -x_{k} e_{ik}$. Since $x_{l} =0$, then \\$\left(V_{0}^{*} \right)_{il} =\left(V_{0} \right)_{il} -\left(x_{k} e_{ik} \right)_{il} =x_{l} -\delta _{kl} x_{k} =0$. It is clear that $\left(V_{0}^{*} \right)_{ik} =W_{ik} =0$ and $V_{0}^{*} +W$ is also a matrix over {\it J}.\\
It is given that $1+V_{0}^{*} c^{2} U\in GL\left(n,R\right).$

Note that in the case $1+J\subseteq R^{*}$ it automatically follows from the equality $V_{0} U=V_{0}^{*} U+x_{k} e_{ll} U$ and inclusion $1+V_0c^2U \in GL(n,R)$ that $1+V_{0}^{*} c^{2} U\in GL\left(n,R\right)$.

As was shown above, matrices $1+Uc^{2} V_{0} $  and  $1+Uc^{2} V_{0}^{*} $ ,  $1+Uc^{2} V_{0}^{*} $  and  $1-Uc^{2} W$ can be expanded into respectively identical constructions of products of transvections and diagonal elements. By taking into account the fact that $V_{0} \equiv V_{0}^{*} \mod J$  and  $V_{0}^{*} \equiv -W\mod J$ we obtain the congruences 
\[\begin{aligned} 1+Uc^{2} V_{0} &\equiv \left(1+Uc^{2} V_{0}^{*}
\right)\mod E\left(n,J\right),    \\  1+Uc^{2} V_{0}^{*} &\equiv
\left(1-Uc^{2} W\right)\mod E\left(n,J\right).\end{aligned}\]
Thus, we have proved the following:
\begin{center}
 $t_{ij} \left(r c^{2}\right)^{g} \in E\left(n,J\right)$, if  $1+V_{0} c^{2} U\in GL\left(n,R\right)$  and $1+\left(V_0-x_ke_{ik}\right)c^2U \in GL(n,R).$
\end{center}
Similarly one can prove that
\begin{center}
 $t_{ij} \left(r c^2\right)^{g'} \in E\left(n,J\right)$, if  $1+V'_{0} c^{2} U'\in GL\left(n,R\right)$  and $1+\left(V'_0-x'_ke_{ik}\right)c^2U' \in GL(n,R).$
\end{center}
In the end, it is proved that 
\begin{center}
 $\left[h,t_{ij} \left( rc^{2}\right)\right]^{g'} =t_{ij} \left(c^{2} r\right)^{g} t_{ij} \left(-c^{2} r\right)^{g'} \subseteq E\left(n,J\right),$ if \end{center} \begin{center}   $1+V_{0} c^{2} U,$ $1+\left(V_0-x_ke_{ik}\right)c^2U,$ $1+\left(V'_0-x'_ke_{ik}\right)c^2U',$ $1+V'_{0} c^{2} U'$
\end{center} 
are contained in $GL\left(n,R\right).$
\par This finishes the proof of 1).
\par Let's note that if in 1) the element $g'$ normalizes the groups $E\left(n,cI\right)$  and  $E\left(n,J\right),$ then
\begin{center}
$\displaystyle \left[h,t_{ij} \left(c^{2} r\right)\right]\subseteq E\left(n,cI\right)\bigcap E\left(n,J\right)$.
\end{center}
Let's prove 2).\\
Let us introduce the following notations.
\begin{center}
 $A=1+Uc^{2} V_{0} $,        $B=(1-Uc^{2} W)^{-1}$,         $C=t_{ij} \left(-c^{2} r\right)$.
\end{center}
Under these notations
\[\begin{aligned}
 \left[g,t_{ij} \left(c^{2} r\right)\right]=t_{ij} \left(c^{2} r\right)^{g} t_{ij} \left(-c^{2} r\right)&=A\left(BC\right)\mbox{and}\\ \left[t_{ij} \left(-c^{2} r\right),g\right]=t_{ij} \left(-c^{2} r\right)t_{ij} \left(c^{2} r\right)^{g}
 &=\left(CA\right)B.
\end{aligned}\]
Obviously, the {\it l}-th column of matrix $Uc^{2} V_{0} $ and the {\it k}-th column of matrix $Uc^{2} W$ are all zero.

Suppose that the equality $k=j=l$ is not satisfied. Then $j\ne k$  or  $j\ne l$.
If $j\ne k$, then the {\it k}-th column of matrix $BC-E$ is all zero. Hence, the commutator $\left[g,t_{ij} \left(c^{2} r\right)\right]$ is a product of two matrices {\it A} and {\it BC}, which satisfy the conditions of Lemma \ref{Petlem5}. In such a case $\left[g,t_{ij} \left(c^{2} r\right)\right]\in \xi GL\left(n,R\right)$  and  $g\in C\left(n,AnnRc^{2} rR\right)$.

Similarly, if $j\ne l$ then the commutator $\left[t_{ij}(-c^2r),g\right]$ is a product of two matrices {\it CA} and {\it B} and the {\it l}-th column of the matrix $CA-E$ is all zero. According to Lemma \ref{Petlem5}, $\left[t_{ij} \left(-c^{2} r\right),g\right]\in \xi GL\left(n,R\right)$ and  $g\in C\left(n,AnnRc^{2} rR\right)$.

Thus, the only case left to consider is $k=j=l$. Then $0=x_{k} =r\left(g^{-1} \right)_{jj} $. Since, according to Lemma \ref{Petlem3}, the diagonal elements of matrix $g\in N$ are not zero divisors, then \\$r=0$,  $AnnRc^{2} rR=Ann0=R$,  $g\in GL\left(n,R\right)\equiv C\left(n,R\right)\equiv C\left(n,AnnRc^{2} rR\right)$.

In the end 2) is proved.

Let's note that the inclusions of Lemma \ref{Petlem6} look like $f\left(c\right)=0$, where {\it f} is a polynomial with coefficients from the ring $R_{n} $  and  $c\in \xi R$. In particular, if 
the inclusions of Lemma \ref{Petlem6} hold for all $c\!\in\! \xi R$ then the elements of $\xi R$ are the roots of polynomial {\it f.}

Let $V_{0} =x_{1} e_{i1} +\cdots +x_{n} e_{in} $, where $x_{l} =0$,  $x_{k} =r\left(g^{-1} \right)_{jk} $,  $1+V_{0} U\in GL\left(n,R\right)$,  $U=ge_{ii} $,  $g\in GL\left(n,R\right)$. Elements of the type $V_0$ are not well defined by matrix $g$. They form a whole class of matrices and elements $r \in R$ act in them as left coefficients of $\left(g^{-1} \right)_{jk} e_{ik}$.

Denote by $R\left(g\right)$ the additive subgroup of the ring {\it R}, generated by all elements $r \in R,$ which appear as left coefficients of the summands $\left(g^{-1} \right)_{jk} e_{ik} $ of matrices $V_{0} =x_{1} e_{i1} +\cdots +x_{n} e_{in} $, where \\$x_{l} =0$,  $x_{k} =r\left(g^{-1} \right)_{jk} $,  $1+V_{0} U\in GL\left(n,R\right),$ $1+(V_0-x_ke_{ik})U \in GL(n,R)$ for some fixed matrix $g\in GL\left(n,R\right)$ and all possible matrices $V_0$ and $1\le k,l\le n$.
\begin{vr}\label{Petozn6} For fixed $1\le i\ne j\le n$, the element $g\in GL\left(n,R\right)$ is called $\left(R,i,j\right)$-stable if $R\left(g\right)=R$.
\end{vr}

In particular, if among $V_{0} =x_{1} e_{i1} +\cdots +x_{n} e_{in}$ there are some that \\$x_{l} =0,$  $x_{k} =r\left(g^{-1} \right)_{jk},$  $V_{0} U\in J\left(R\right)e_{ii},$ $x_kg_{ki} \in J(R),$ then $RV_{0} U\in J\left(R\right)e_{ii}$,  $1+RV_{0} U\subseteq 1+J\left(R\right)e_{ii} \subseteq GL\left(n,R\right),$ $1+R(V_0-x_ke_{ik})U \in GL(n,R)$ and, henceforth, $R\left(g\right)$ contains the left ideal {\it Rr}. If at the same time $r\in R^{*},$ then $R\left(g\right)=R$ and {\it g} is a $\left(R,i,j\right)$-stable element.

Therefore, the matrix $g\in GL\left(n,R\right)$, for which the inclusion 
\begin{center}
$\left(g^{-1} \right)_{jk} g_{ki} \in J\left(R\right)$
\end{center}
holds, is $\left(R,i,j\right)$-stable. In order to show this, it is enough to choose $V_{0} \in R\left(g^{-1} \right)_{jk}e_{ik} $ and use the inclusion $1+V_{0} U\subseteq GL\left(n,R\right).$
\begin{vr}\label{Petozn7} Element $g\in GL\left(n,R\right)$ which is $\left(R,i,j\right)$-stable for all \\$1\le i\ne j\le n$ such that $E\left(n,R\right)=\left\langle t_{ij} \left(R\right)\right\rangle $ is called {\it R}-stable.
\end{vr}

The example of {\it R}-stable element is the matrix {\it g} with $g_{ij} =g_{ji} =0$ for all $1\le i\ne j\le n$, where {\it j} is fixed. As in this case {\it g} is $\left(R,i,j\right)$-stable and $\left(R,j,i\right)$-stable element.

Since $E\left(n,R\right)=\left\langle t_{ij} \left(R\right),t_{ji} \left(R\right)\left|{\rm \; }1\le i\ne j\le n,{\rm \; }j-\mbox{fixed}\right. \right\rangle$, then {\it g} is {\it R}-stable element.

It is not hard to see that $\left(R,i,j\right)$-stability and, as a consequence, {\it R}-stability of matrices is preserved when factoring the ring {\it R}.
\par Let $I$ be an arbitrary ideal of the ring $R,$ and $N$ is a subgroup of $GL\left(n,R\right)$ invariant with respect to $E\left(n,R\right)$, which does not contain non-identity transvections.
\begin{vr}\label{Petozn8} Element $g\in GL\left(n,R\right)$, for which from $g\in C_{I}$ it follows that $\left[g,E\left(n,J\right)\right]\subseteq E\left(n,I\right)\bigcap E\left(n,J\right),$ and from $g\in N,$ follows that $g\in \xi GL\left(n,R\right)$ is called stable.
\end{vr}

As shall be shown below, for instance, elements of the group $C\left(n,J\left(R\right)\right)$ are stable.
\begin{vr} We shall say that the class of invertible matrices $L_{g},$ which contains an identity, up to transvections, approximates the element $g \in GL(n,R)$, if for $g\in C_{I}$ there exists an element $g_{I} \in L_{g}$ such that \\$g_{I} \in E\left(n,I\right)g^{E\left(n,R\right)} E\left(n,I\right),$ and for $g\in N$ there exists an element $g_{N} \in L_{g}$ such that 
$g_{N} \in \left[g,t_{ij} \left(R^{*} \right)\right]^{E\left(n,R\right)}$ for some $1\le i\ne j\le n$.
\end{vr}
\begin{lr}\label{Petlem8} From the stability of elements $L_{g}$ one has the stability of element {\it g}.
\end{lr}

{\it\textbf{Proof.}} If $g\in C_{I}$, then there exists a stable element \\$g_{I} \in E\left(n,I\right)g^{E\left(n,R\right)} E\left(n,I\right)\subseteq C_{I}$. From stability of element $g_{I} $ it follows that the following inclusions hold
\begin{center}
 $\left[g_{I} ,E\left(n,J\right)\right]\subseteq E\left(n,I\right)\bigcap E\left(n,J\right)$   and  $\left[g,E\left(n,J\right)\right]\subseteq E\left(n,I\right)\bigcap E\left(n,J\right)$
\end{center}
for all ideals {\it I}{\it , }{\it J} of the ring {\it R}.

Similarly, if $g\in N$, then there exists a stable element $g_{N}$ such that \\$ g_{N}\in \left[g,t_{ij} \left(R^{*} \right)\right]^{E\left(n,R\right)} \subseteq N$. From stability of element $g_{N}$ it follows that $g_{N} \in \xi GL\left(n,R\right)$. Hence, there exists an element $r\in R^{*}$ such that $\left[g,t_{ij} \left(r\right)\right]\in \xi GL\left(n,R\right)$. According to Corollary \ref{Petnas3} we have \[g\in C\left(n,AnnRrR\right)=C\left(n,AnnR\right)=\xi GL\left(n,R\right).\]

\begin{theo}\label{Petthem1} Let {\it R} be an associative ring with identity. Elements of the group $GL\left(n,\right.$ $\left.R\right),$ which, up to transvections, are approximated by the class of {\it R}-stable matrices, are stable. If, up to transvections, all elements of the group $GL\left(n,R\right)$ are approximated by classes of {\it R}-stable matrices, then {\it R} is a stable ring.
\end{theo}

{\it\textbf{Proof.}} According to Lemma \ref{Petlem8} it is sufficient to prove that {\it R}-stable elements are in fact stable.
\par Let $g\in GL\left(n,R\right)$ be an {\it R}-stable element. Then {\it g} is a $\left(R,i,j\right)$-stable element for all pairs $1\le i\ne j\le n$ such that $E\left(n,R\right)=\left\langle t_{ij} \left(R\right)\right\rangle$. Let's fix one such pair $1\le i\ne j\le n$. Then the additive group 
\[\begin{aligned}
R\left(g\right)=\langle r\|&1\!+\!V_{0} U\in GL\left(n,R\right),\! 1+(V_0-x_ke_{ik})U \in GL(n,R), \\&V_{0}\!=\!x_{1} e_{i1}\! +\!\cdots\!+x_{n} e_{in} ,x_{l}\!=\!0,x_{k}\!=\!r\left(g^{-1} \right)_{jk} \rangle _{Z}\!=R.
\end{aligned}\]

\par According to this, $1+\alpha \in R^{*}$ and $1+\alpha-x_k g_{ki} \in R^{*}.$
\par Let {\it I}{\it , }{\it J} be ideals of {\it R}, $g\in C_{I} $, $r\in J$,  $U=ge_{ii},$ $V=e_{ij}g^{-1},$ $\alpha e_{ii} =V_{0} U$. Obviously $x_{i}-\alpha \in I.$
We shall use Lemma \ref{Petlem7}.

In the case when the equality $l=k=j$ fails, we put $c=1$ ,  $g'=1,$ $x'_{i} =(1-\delta_{ik})\left(1-\delta _{il} \right)\alpha,$ where $x'_{k}=r\delta _{jk},$ $x'_{s}=x_{s} $ for all $1\le s\ne i,k\le n.$
\par Then $x'_{l}=0.$ This is obvious if $l\neq i,k$ or $l=i.$ If $l=k,$ then \\$x'_l=x'_k=r\delta_{jl}=0$ when $l=k\neq j.$ Moreover,  $x'_t-x_t \in I $ for all $1 \leq t \neq i \leq n,$ $x'_i-x_i=\alpha-x_i-\delta_{il}\alpha \in I.$ Since $U'=e_{ii},$ we have
\begin{center}$V'_0U'=x'_ie_{ii},$
$(V'_0-x'_ke_{ik})U'=V'_0 U'-r\delta_{jk} e_{ik} e_{ii}=V'_0
U'-r\delta_{jk}\delta_{ik}e_{ii}=V'_0 U',$ 
\end{center}
and $1+(V'_0-x'_ke_{ik})U'=1+V'_0U'=1+x'_ie_{ii} \in R^{*},$ as $x'_i=0$ or $x'_i=\alpha.$
\par In view of Lemma 7,
\begin{center} $\left[g,t_{ij} \left(r\right)\right]\subseteq
E\left(n,I\right)\bigcap E\left(n,J\right)$ 
\end{center}  
and, as a consequence,
\begin{center}$\left[g,t_{ij}
\left(R\left(g\right)\right)\right]\subseteq
E\left(n,I\right)\bigcap E\left(n,J\right).$ \end{center}

If $l=k=j$, then $0=x_{l} =x_{k} =r\left(g^{-1} \right)_{jk} =r\left(g^{-1} \right)_{jj}$ and, as a consequence, $r\in I.$ It is not hard to see that $V_{0} =rV$ satisfies all conditions, defined by the additive group $R(g).$ Recall, that in Lemma \ref{Petlem6}
\begin{center} $a=Ue_{il} -U_{l} e_{ll} $ , $b=e_{li} V_{0} =re_{li} V$.
\end{center}

Taking into account that $VU=0$  and  $\alpha e_{ii} =V_{0} Uб$ we obtain the equalities \[\alpha =0, ba=-re_{li} Vu_{l} e_{ll} =-re_{jj} Vu_{l} e_{jj}=-r\left(g^{-1} \right)_{jj} u_{l} e_{jj} =0.\]
Henceforth, $\gamma =\left(1+ab\right)^{-1} =1-ab$ ,  $\left(1-\gamma \right)a=b\left(1-\gamma \right)=0,$  $d_{l} =d_{k} =1,$
\begin{center}
$\displaystyle t_{ij} \left(r\right)^{g}=T\left(g\right)=\left[1-b,1+a\right]\left(1+u_{l} b\right)$.
\end{center}

If $r\in I\bigcap J$ then the matrices $1-b$ ,  $1+u_{l} b$ are contained in the group $E\left(n,I\right)\bigcap$ $ E\left(n,J\right).$
\par Thus, in the case $l=k=j,$ we have the inclusion
\begin{center}
$\left[g,t_{ij} \left(r\right)\right]\subseteq E\left(n,I\right)\bigcap E\left(n,J\right).
$ \end{center}
as well.

Therefore, it is proven that in all cases 
\begin{center}$\left[g,t_{ij} \left(J\bigcap R\left(g\right)\right)\right]\subseteq E\left(n,I\right)\bigcap E\left(n,J\right)$.
\end{center}
The $\left(R,i,j\right)$-stability of element {\it g} implies
\begin{center}
 $\left[g,t_{ij} \left(J\right)\right]\subseteq E\left(n,I\right)\bigcap E\left(n,J\right)$     and     $\left[g,t_{ij} \left(R\right)\right]\subseteq E\left(n,I\right)$.
\end{center}
Hence, the {\it R}-stability of element {\it g} implies the inclusion 
$$\left[g,E\left(n,R\right)\right]\subseteq E\left(n,I\right).$$

Taking into account the fact that 
\begin{center}
 $E\left(n,J\right)=\left\langle t_{ij} \left(J\right)^{E\left(n,R\right)} \right\rangle $       and      $\left[E\left(n,I\right),E\left(n,J\right)\right]\subseteq E\left(n,I\right)\bigcap E\left(n,J\right)$
\end{center}
we obtain the necessary inclusion
\begin{center}
$\displaystyle \left[g,E\left(n,J\right)\right]\subseteq
E\left(n,I\right)\bigcap E\left(n,J\right)$.
\end{center}
Let {\it N} be a group, invariant with respect to $E\left(n,R\right)$, not containing non-identity transvections. If $g\in N$, in accordance with Lemma \ref{Petlem7}, we have \\$g\in C\left(n,AnnRrR\right)$ for all generators {\it r} of the additive group $R\left(g\right)=R.$ Hence, $g\in C\left(n,AnnR\right)=\xi GL\left(n,R\right)$. Thus, it is proved that {\it R}-stability of elements of the group $GL\left(n,R\right)$ implies their stability.

If, up to transvections, all elements of the group $GL\left(n,R\right)$ are approximated by a class of stable matrices, then
\begin{center}
 $\left[C\left(n,I\right),E\left(n,R\right),E\left(n,J\right)\right]\subseteq E\left(n,I\right)\bigcap E\left(n,J\right)$  and  $N\in \xi GL\left(n,R\right)$.
\end{center}
If one puts $J=R$, then it follows that {\it R} is a weakly-commutator ring. Moreover, {\it R} is a partially normal ring. Since {\it R}-stability is preserved when taking quotients, then all quotient rings of {\it R} are partially normal as well. According to Lemma \ref{Petlem2}, {\it R} is a stable ring.

Theorem \ref{Petthem1} implies that {\it commutative rings with identity are stable} \cite{Pet6, Pet7, Pet10, Golubchik}.

Indeed, let {\it R} be a commutative ring with identity, $g\in GL\left(n,R\right)$. For arbitrary $1\le i\ne j\le n$ define 
\begin{center}
$\displaystyle V_{0} =\left(g^{-1} \right)_{jk} \left(g_{ti} e_{ik}
-g_{ki} e_{it} \right)=g_{ti} \left(g^{-1} \right)_{jk} e_{ik}
-g_{ni} \left(g^{-1} \right)_{jk} e_{it} $,
\end{center}
where $1\le k\ne t\le n$.

Then $V_{0} U=0$ and $Rg_{ti} \subseteq R\left(g\right)$.
\par By interchanging {\it k} and {\it t} we have $R\left(g\right)=R$. This means that all elements of the group $GL\left(n,R\right)$ are $\left(R,i,j\right)$-stable and, as a consequence, {\it R}-stable. According to Theorem 1, {\it R} is a stable ring.

\begin{theo}\label{Petthem2} Let {\it R} be an associative ring with identity, $g\in GL\left(n,R\right)$ and at least one element of the matrix {\it g} belongs to the radical $J\left(R\right)$. Then {\it g} is a stable element.
\end{theo}

{\it\textbf{Proof.}} Let $g=\left(g_{ij} \right)\in GL\left(n,R\right)$ and $g_{ij} \in J\left(R\right).$
We introduce the notation $g_{1} =e_{1} ege_{2},$ where
\begin{center}
 $e=t_{i1} \left(\left(g^{-1} \right)_{j1} \right)\cdots t_{in} \left(\left(g^{-1} \right)_{jn} \right),$
 $\alpha=-\left(1+g_{ij}-(g^{-1})_{ji}g_{ij}\right)^{-1},$     $e_{1} =t_{1i} \left(g_{1j}^{} \alpha \right)\cdots t_{ni} \left(g_{nj}^{} \alpha \right)$,  $e_{2} =t_{j1} \left(\alpha \left(1-\left(g^{-1} \right)_{ji} \right)g_{i1} \right)\cdots t_{jn} \left(\alpha \left(1-\left(g^{-1} \right)_{ji} \right)g_{in} \right).$
\end{center}

It is not hard to see that 
\begin{center}
 $\left(g_{1} \right)_{il} =\left(g_{1} \right)_{sj} =0$    for all  $1\le l\ne j\le n$ ,   $1\le s\ne i\le n$ .
\end{center}
Let {\it I} be an ideal of {\it R}, $g\in C_{I} $.

If $i=j$ then $e,e_{1} ,e_{2} $ are contained in $E_{I} $,  $g_{1} \in C_{I} $ and $g_{1} $ satisfies the condition of {\it R}-stability.

If $i\ne j$, then there exists $e_{0} =t_{ji} \left(1\right)t_{ij} \left(-1\right)t_{ji} \left(1\right)$ such that $e_{0} e_{1} ee_{2} \in E\left(n,I\right)$,  $e_{0} g_{1} \in C_{I}$ and satisfies the condition of {\it R}-stability. 

Thus, it is proved that $g\in C_{I}$ and is approximated, up to transvections, by {\it R}-stable matrices. Therefore, {\it g} satisfies the condition of {\it R}-stability.

If $g\in N$ then $1+g_{ki}^{-1} g_{ij} \in R^{*} $ and, according to Lemma \ref{Petlem7}, $g\in C\left(n,AnnR\right)=\xi GL\left(n,R\right)$. 
Thus, it is proved that {\it g} is a stable element.

In particular the following inclusions hold
\[\begin{aligned}
      \left[C\left(n,J\left(R\right)\right)\bigcap C\left(n,I\right),E\left(n,J\right)\right]&\subseteq E\left(n,I\right)\bigcap E\left(n,J\right),  \\
 N\bigcap C\left(n,J\left(R\right)\right)&\subseteq \xi
      GL\left(n,R\right).
\end{aligned}\]

If we put $I=R$, then $\left[C\left(n,J\left(R\right)\right),E\left(n,J\right)\right]\subseteq E\left(n,J\right)$.
\begin{lr}\label{Petlem9} Let {\it R} be an associative ring with identity, $J\left(R\right)$ - radical, ${\raise0.7ex\hbox{$ R $}\!\mathord{\left/{\vphantom{R J\left(R\right)}}\right.\kern-\nulldelimiterspace}\!\lower0.7ex\hbox{$ J\left(R\right) $}}$ - partially normal ring. Then {\it R} is a partially normal ring.
\end{lr}

{\it\textbf{Proof.}} Let {\it N} be a subgroup of the group $GL\left(n,R\right)$, invariant with respect to $E\left(n,R\right)$ and not containing non-identity transvections. If $\Lambda _{J\left(R\right)} N$ contains non-identity transvection $\Lambda _{J\left(R\right)} t_{ij} \left(r\right)$,  $r\notin J\left(R\right)$, then $t_{ij} \left(r\right)\in hC_{J\left(R\right)}$ where $h\in N$. Thus, $h\in t_{ij} \left(r\right)C_{J\left(R\right)}$ and at least one element of the matrix {\it h} is contained in $J\left(R\right)$. According to Theorem \ref{Petthem2}, $h\in \xi GL\left(n,R\right)$ and $r\in J\left(R\right),$ contradicting the assumption. Therefore, the group $\Lambda _{J\left(R\right)} N$ does not contain non-identity transvections. Since ${\raise0.7ex\hbox{$ R $}\!\mathord{\left/{\vphantom{R J\left(R\right)}}\right.\kern-\nulldelimiterspace}\!\lower0.7ex\hbox{$ J\left(R\right) $}} $ is a partially normal ring, then $\Lambda _{J\left(R\right)} N\in \xi GL\left(n,{\raise0.7ex\hbox{$ R $}\!\mathord{\left/{\vphantom{R J\left(R\right)}}\right.\kern-\nulldelimiterspace}\!\lower0.7ex\hbox{$ J\left(R\right) $}} \right)$. This means that $N\subseteq C\left(n,J\left(R\right)\right)$. By Theorem \ref{Petthem1}, $N\subseteq \xi GL\left(n,R\right).$ Thus, it is proved that {\it R} is a partially normal ring.

\begin{lr}\label{Petlem10} Let {\it R} be an associative ring with identity, $J\left(R\right)$- radical, ${\raise0.7ex\hbox{$ R $}\!\mathord{\left/{\vphantom{R  J\left(R\right)}}\right.\kern-\nulldelimiterspace}\!\lower0.7ex\hbox{$ J\left(R\right) $}} $-normal ring. Then all the quotient rings of the ring {\it R} are partially normal rings.
\end{lr}

{\it\textbf{Proof.}} Let $\overline{R}$ be some quotient ring of the ring {\it R}. Since under the epimorphism of rings, preimages of maximal one-sided ideals are maximal one-sided ideals, then $\overline{J\left(R\right)}\subseteq J\left(\overline{R}\right).$ Hence, the epimorhism of rings $R \rightarrow \overline{R} \rightarrow
{\raise0.7ex\hbox{$ \overline{R} $}\!\mathord{\left/{\vphantom{R
J\left(R\right)}}\right.\kern-\nulldelimiterspace}\!\lower0.7ex\hbox{$ J\left(\overline{R}\right) $}} $ induces the epimorhism ${\raise0.7ex\hbox{$ R $}\!\mathord{\left/{\vphantom{R
J\left(R\right)}}\right.\kern-\nulldelimiterspace}\!\lower0.7ex\hbox{$ J\left(R\right) $}}  \rightarrow {\raise0.7ex\hbox{$ \overline{R} $}\!\mathord{\left/{\vphantom{R
J\left(R\right)}}\right.\kern-\nulldelimiterspace}\!\lower0.7ex\hbox{$J\left(\overline{R}\right) $}}$ and, as a consequence,
${\raise0.7ex\hbox{$ \overline{R}$}\!\mathord{\left/{\vphantom{\overline{R}
J\left(\overline{R}\right)}}\right.\kern-\nulldelimiterspace}\!\lower0.7ex\hbox{$ J\left(\overline{R}\right) $}} $  is a quotient ring of the normal ring ${\raise0.7ex\hbox{$ R
$}\!\mathord{\left/{\vphantom{R J\left(R\right)}}\right.\kern-\nulldelimiterspace}\!\lower0.7ex\hbox{$ J\left(R\right) $}}.$
Since all quotients of normal ring are partially normal rings, then ${\raise0.7ex\hbox{$
\overline{R} $}\!\mathord{\left/{\vphantom{\overline{R} 
J\left(\overline{R}\right)}}\right.\kern-\nulldelimiterspace}\!\lower0.7ex\hbox{$ J\left(\overline{R}\right) $}} $ is a partially normal ring.
In view of Lemma \ref{Petlem9}, $\overline{R}$ is a partially normal ring.

From Theorem \ref{Petthem1} one has
\begin{sr}\label{Petnas5} \cite{Pet27} Let {\it R} be an associative ring with identity, $J\left(R\right)$ - radical, ${\raise0.7ex\hbox{$ R $}\!\mathord{\left/{\vphantom{R J\left(R\right)}}\right.\kern-\nulldelimiterspace}\!\lower0.7ex\hbox{$ J\left(R\right) $}} $ - a stable ring. Then {\it R} is a stable ring.
\end{sr}

{\it\textbf{Proof.}} Since  ${\raise0.7ex\hbox{$ R $}\!\mathord{\left/{\vphantom{R J\left(R\right)}}\right.\kern-\nulldelimiterspace}\!\lower0.7ex\hbox{$ J\left(R\right) $}} $ is a stable ring, then ${\raise0.7ex\hbox{$ R $}\!\mathord{\left/{\vphantom{R J\left(R\right)}}\right.\kern-\nulldelimiterspace}\!\lower0.7ex\hbox{$ J\left(R\right) $}} $ is commutator and normal ring.

From commutatorness of the ring  ${\raise0.7ex\hbox{$ R $}\!\mathord{\left/{\vphantom{R J\left(R\right)}}\right.\kern-\nulldelimiterspace}\!\lower0.7ex\hbox{$ J\left(R\right) $}} $ it follows that for arbitrary ideals {\it I} and {\it J} of the ring {\it R} one has the inclusions
$$\left[C\left(n,I\right),E\left(n,J\right)\right]\subseteq \left(E\left(n,I\right)\bigcap E\left(n,J\right)\right)C_{I\bigcap J\left(R\right)}. $$
\par According to Theorem \ref{Petthem2} $$\left[C_{I\bigcap J\left(R\right)} ,E\left(n,J\right)\right]\subseteq E\left(n,I\right)\bigcap E\left(n,J\right).$$
\par     Therefore,
\[\begin{aligned}
\left[C\left(n,I\right),E\left(n,J\right),E\left(n,J\right)\right] &\subseteq \left[\left(E\left(n,I\right)\bigcap
E\left(n,J\right)\right)C_{I\bigcap J\left(R\right)} ,E\left(n,J\right)\right] \subseteq \\ &\subseteq E\left(n,I\right)\bigcap
E\left(n,J\right).
\end{aligned}\]
\par In particular, if $J=R$, then $$\left[C\left(n,I\right),E\left(n,R\right),E\left(n,R\right)\right]\subseteq E\left(n,R\right).$$
\par This means that {\it R} is a weakly-commutator ring.

Since ${\raise0.7ex\hbox{$ R $}\!\mathord{\left/{\vphantom{R J\left(R\right)}}\right.\kern-\nulldelimiterspace}\!\lower0.7ex\hbox{$ J\left(R\right) $}} $ is a normal ring, then, in view of Lemma \ref{Petlem9}, all the quotient rings of the ring {\it R} are partially normal.

According to Lemma \ref{Petlem2}, weakly-commutator rings all quotients of which are partially normal, are stable.

\begin{vr} Vector  $\left(r_{1} ,\ldots ,r_{n} \right)$ is called unimodular in $R^{n} $, if there are elements $t_{1} ,\ldots ,t_{n} $ in the ring {\it R} such that $t_{1} r_{1} +\cdots +t_{n} r_{n} =1.$
\end{vr}

Obviously, in the case $n=1$ unimodular vectors are exactly the left invertible elements of the ring {\it R}.

From the inclusion $1+J(R) \subseteq R^*$ one has that vectors, unimodular modulo the radical $J(R)$, are in fact unimodular.
\begin{vr} Let $n\ge 2$. The associative ring {\it R} is said to satisfy the condition of stability of rank $n-1$, if for an arbitrary unimodular vector $\left(r_{1} ,\ldots ,r_{n} \right)$ there are elements $s_{2} ,\ldots ,s_{n} $ in {\it R} such that the vector $\left(r_{2} +s_{2} r_{1} ,\ldots ,r_{n} +s_{n} r_{1} \right)$ is unimodular in $R^{n-1} $.
\end{vr}

It is not hard to see that the direct sum of rings, which satisfy the condition of stability of rank $n-1$, satisfies the condition of stability of rank $n-1$ as well. Analogously, rings which modulo the radical satisfy the condition of stability of rank $n-1$, satisfy the condition of stability of rank $n-1$.

It is known \cite{Pet4}, that if {\it R} satisfies the condition of stability of rank $m-1,$ then {\it R} satisfies the condition of stability of rank $n-1$ for all $n\ge m$.

Let us consider associative rings $R$ with identity such that for an arbitrary element $e \in R$ there exist elements $r_1, r_2 \in R^*$ such that $r_1er_2$ is an idempotent of the ring $R$.

In particular, matrix rings over skew fields satisfy the aforementioned condition. Indeed, by means of elementary transformations an arbitrary matrix in the ring of matrices over skew field can be brought to the idempotent \\ $diag\left(1,\ldots ,1,0,\ldots ,0\right)$.

Let us now prove that the considered rings satisfy the condition of stability of rank 1. Let $(r,e)$ be a unimodular element. Find an element $s \in R$ such that $e+sr \in R^*$. At first we consider the case when $e$ is an idempotent of $R$.

Let $1=\alpha r+\beta e$, where $\alpha$ and $\beta$ are some elements of the ring $R$, $e^2=e$ is an idempotent of the ring $R$. Since $\left(1-e\right)\beta e$ is a nilpotent element then by putting $s=\left(1-e\right)\alpha $ we obtain 
\[e+sr=e+\left(1-e\right)\alpha r=e+\left(1-e\right)\left(1-\beta e\right)=1-\left(1-e\right)\beta e\in R^{*}.\]
If $e$ is an arbitrary element of the ring $R$ such that $(r,e)$ is unimodular and element $r_1er_2$ is an idempotent of the ring $R$, where $r_1 , r_2 \in R^*$, then $(r_1rr_2,r_1er_2)$ is unimodular and, according to the argument above, there exists an element $s_1 \in R$ such that $r_1er_2 +s_1r_1rr_2 \in R^*$. By taking $s=r_{1}^{-1}s_1r_1$ we obtain $e+sr \in R^*$.

In the end, it is proven that the associative rings with identity which modulo the radical are the direct sum of rings with the property that their elements, multiplied by invertible elements of the ring, can be turned into idempotents, satisfy the condition of stability of rank 1.

\begin{vr} Associative ring {\it R} with identity is called semilocal if ${\raise0.7ex\hbox{$R$}\!\mathord{\left/{\vphantom{R J\left(R\right)}}\right.\kern-\nulldelimiterspace}\!\lower0.7ex\hbox{$ J\left(R\right) $}} $ is a direct sum of full matrix rings over skew fields.
\end{vr}

Semilocal ring is a classic example of the ring which satisfies the condition of stability of rank 1.

%
%

{\it Associative rings with identity which satisfy the condition of stability of rank $n-1>1$ are stable} \cite{Pet1}.

Indeed, let $g\in GL\left(n,R\right)$. Vector $\left(g_{1n} ,g_{2n} ,\ldots ,g_{nn} \right)$ is obviously unimodular. Therefore, there are elements $k_{2},\ldots ,k_{n} $ such that 
\[\left(g_{2n} +k_{2} g_{1n} ,\ldots ,g_{nn} +k_{n} g_{1n} \right)\]
 is also a unimodular vector. Hence, there are elements $s_{2},\ldots ,s_{n} $ in $g_{1n} R$ such that 
\[g_{1n} +s_{2} \left(g_{2n} +k_{2} g_{1n} \right)+\cdots +s_{n} \left(g_{nn} +k_{n} g_{1n} \right)=0.\]
\par Let $e_{1} =t_{21} \left(k_{2} \right)\cdots t_{n1} \left(k_{n} \right)$,  $e_{2} =t_{12} \left(s_{2} \right)\cdots
t_{1n} \left(s_{n} \right)$, $g_{1} =e_{2} g^{e_{1} }$, \\$g_{2}=\left[g^{e_{1} },t_{n2} \left(1\right)\right]^{e_{2} } =t_{n2}\left(1\right)^{g_{1} } t_{n2} \left(-1\right)^{e_{2}}.$
\par Then $(g_1)_{1n}=(g_2)_{1n}=0.$ By Theorem \ref{Petthem2}, elements  $g_1$ and $g_2$ are stable.
Obviously, $e_2 \in E_I$ if $g \in C_I$ and $g_2 \in N$ if $g \in N.$ Since elements $g_1$ and $g_2,$ up to transvections, approximate the matrix $g,$ then, according to Lemma \ref{Petlem8}, $g$ is a stable element. Thus, it is proved that $R$ is a stable ring.
\par We shall need the useful
\begin{lr}\label{Petlem11} Let {\it R} be an associative ring with identity, {\it I}, {\it J} - ideals of the ring {\it R}, {\it N} - subgroup of $GL\left(n,R\right)$ invariant with respect to $E\left(n,R\right)$ and not containing non-identity transvections, $g\in GL\left(n,R\right)$ and there exists an element $e\in R$ such that $e^{2} -e\in J\left(R\right)$,  $g_{jk} \in eR$,  $e\in g_{jk} R$ for some $1\le k,j\le n.$ Then $\left[g^{-1} ,t_{ij} \left(J\right)\right]\subseteq E\left(n,I\right)\bigcap E\left(n,J\right)$ if $g\in C_{I}$, $1\le i\ne k,j\le n$  and $g\in \xi GL\left(n,R\right)$ if $g\in N$ and  $j=k$.
\end{lr}

{\it\textbf{Proof.}} Let $e=g_{jk} r$,  $g_{e} =t_{ki}\left(rg_{ji} \right)g^{-1},$ where $r\in R$,  $1\le i\ne k,j\le n$.
It is not hard to see that $g_{e}^{-1} =gt_{ki} \left(-rg_{ji} \right)$ and
\begin{center}
 $\left(g_{e}^{-1} \right)_{jk} =g_{jk} \in eR,$ \quad $\left(g_{e}^{-1} \right)_{ji} =g_{ji} -g_{jk} rg_{ji} =\left(1-g_{jk} r\right)g_{ji} \in \left(1-e\right)R$.
\end{center}
\par Since $\left(1-e\right)\left(g_{e}^{-1} \right)_{jk} \left(g_{e}\right)_{ki} \in J\left(R\right)$  and  $e\left(g_{e}^{-1}
\right)_{ji} \left(g_{e} \right)_{ii} \in J\left(R\right),$
 then 
\begin{center} $R\left(1-e\right)\subseteq R\left(g_{e} \right)$ and
$Re\subseteq R\left(g_{e} \right).$ 
\end{center}
Therefore, $R\subseteq Re+R\left(1-e\right)\subseteq R\left(g_{e}\right)\subseteq R$. Hence, $R\left(g_{e} \right)=R$ and $g_{e} $ is an $\left(R,i,j\right)$-stable element for all $1\le i\ne k,j\le n$.

\par If $g \in C_I,$ then $g_{ji} \in I,$ $g_{e}\in C_I,$ $\left[g_{e}, t_{ij}(J)\right]\subseteq E(n,I) \bigcap E(n, J)$ and, as a consequence, 
\begin{center} $\left[g^{-1},t_{ij}(J)\right]\subseteq E(n,I) \bigcap E(n, J).$
\end{center}
\par If $g \in N$ and $j=k$, then as was shown in Lemma \ref{Petlem9}, $\Lambda_{J(R)}N$ does not contain non-identity transvections. In accordance with Lemma \ref{Petlem3}, diagonal elements of matrices $\Lambda_{J(R)}N$ are zero divisors free. Hence, $1-e \in J(R),$ $e \in R^{*}$ and $g_{jj} \in R^{*}.$ According to Corollary \ref{Petnas4}, $g \in \xi GL(n, R).$
\par Let's note that if Lemma \ref{Petlem11} holds for all $1\le i\ne k,j\le n$, then $\left[g^{-1} ,E_{J} \right]\subseteq
E\left(n,I\right)\bigcap E\left(n,J\right)$ and, as a consequence, when $J=R$
\[\left[g^{-1} ,E\left(n,R\right)\right]\subseteq E\left(n,I\right)\bigcap E\left(n,J\right).\]
Thus, 
\[\left[g^{-1} ,E\left(n,J\right)\right]\subseteq E\left(n,I\right)\bigcap E\left(n,J\right), \;\mbox{if}\; g\in C_{I} \;\mbox{and}\; g\in \xi GL\left(n,R\right), \;\mbox{if}\; g \in N.\] 

This means that {\it g} is a stable element.
\begin{vr} Associative ring $ R$ with identity is called von Neumann regular if for an arbitrary element $a\in R$ there exists an element $a'\in R$ such that $aa'a=a$.
\end{vr}

It turns out \cite{Pet12, Pet17} that {\it von Neumann regular rings are stable}.

Indeed, let {\it R} be a von Neumann regular ring and $g\in GL\left(n, R\right).$ Let $a=g_{jk} $ for arbitrary $1\le k,j\le n$. Then there exists an element $a'\in R$ such that $aa'a=a$. Let $e=aa'$. Obviously, $ea=a,$ $e^{2} =eaa'=e,$   $g_{jk} =a\in eR$ and $e\in g_{jk} R$. According to Lemma \ref{Petlem11},
\begin{center}  
$\left[g^{-1} ,t_{ij} \left(J\right)\right]\subseteq E\left(n,I\right)\bigcap E\left(n,J\right)$ if  $g\in C_{I} $ and $g\in \xi GL\left(n,R\right)$ if $g \in N,$
\end{center}
where {\it I}, {\it J} - ideals of the ring {\it R} and {\it N} is a subgroup of $GL\left(n,R\right)$ invariant with respect to $E\left(n,R\right)$ and does not contain non-identity transvections.

\par Thus, {\it g} is a stable element and, as a consequence, {\it R} is a stable ring.
\par Let us present an author's
\begin{vr} Associative ring $R$ with identity is called nearly local if for an arbitrary $a\in R$ there exists an element $a'\in R$ such that \\$\left(1+a'a\right)\left(1-a'+aa'\right)=0.$
\end{vr}

Obviously, local rings with identity and their direct and Cartesian products are nearly local rings.

\begin{theo}\label{Petthem3} Nearly local rings are stable.
\end{theo}

{\it\textbf{Proof.}}
     Let  $g\in GL\left(n,R\right)$, $r\in R$, $a=\left(g^{-1} \right)_{ii} g_{ii}$ where $1\le i\le n.$
Define
\begin{center} $g_{0} =t_{1i} \left(g_{1i} r\left(g^{-1} \right)_{ii} \right)\cdots t_{ni} \left(g_{ni} r\left(g^{-1} \right)_{ii} \right)$, $g_{1} =g^{-1} g_{0}^{-1} $, $g_{2} =\left(g^{-1} \right)^{g_{0} }.$
\end{center}

It is not hard to see that 
\begin{center}
$\left(g_{1} \right)_{ii} =\left(1-r+ar\right)\left(g^{-1}\right)_{ii} =\left(g_{2} \right)_{ii} .$
\end{center}

Let $a'$ be an element of $R,$ for which $\left(1+a'a\right)\left(1-a'+aa'\right)=0.$

Define $e=\left(1-a'+aa'\right)a$. Then $1-e=\left(1-a\right)\left(1+a'a\right)$.
Obviously, \[\left(1+a'a\right)e=0, \quad 
\left(1-e\right)\left(1-a'+aa'\right)=0 \quad \mbox{and} \quad \left(1-e\right)e=0.\]

This means that $e^{2} =e$  and  $\left(1-e\right)^{2} =1-e$.
\par Therefore, $eR=\left\{t \in R \left|(1-e)t=0\right.\right\}$ and $(1-e)R=\left\{t \in R\left|et=0 \right. \right\}.$
\par Hence, $1-a'+aa'\in eR$. Suppose that in the definition of elements $g_{0} ,g_{1}$ and $g_2$ element $r=a'$. Then
\begin{center}
$\left(g_{1} \right)_{ii} =\left(g_{2} \right)_{ii}=\left(1-a'+aa'\right)\left(g^{-1} \right)_{ii} \in eR.$
\end{center}
If $a=\left(g^{-1} \right)_{ii} g_{ii} $, then 
\[e=\left(1-a'+aa'\right)a=\left(1-a'+aa'\right)\left(g^{-1}\right)_{ii} g_{ii} =\left(g_{1} \right)_{ii} g_{ii} =\left(g_{2}\right)_{ii} g_{ii}.\]
Thus, Lemma \ref{Petlem11} can be applied to elements $g_{1} $ and $g_{2}$.

\par In particular, if $g\in C_{I} $, then $g_{0} \in C_{I} $, $g_{1} \in C_{I}, $
\begin{center}
$\left[g_{1}^{-1} ,E\left(n,J\right)\right]\subseteq E\left(n,I\right)\bigcap E\left(n,J\right)$ and
$\left[g,E\left(n,J\right)\right]\subseteq E\left(n,I\right)\bigcap E\left(n,J\right),$
\end{center}
where \textit{I, J} -- ideals of \textit{R}.

\par If $g \in N,$ then $g_{2} \in N$ and $g_{2} \in \xi GL\left(n,R\right)$, $g\in \xi GL\left(n,R\right)$, where \textit{N} is a subgroup of $GL\left(n,R\right)$ invariant with respect to $E\left(n,R\right)$ and does not contain non-identity transvections.

Thus, it is proved that all elements of the group $GL\left(n,R\right)$ are stable. According to Theorem \ref{Petthem1},  \textit{R} is a stable ring.

Similarly, one can prove that the {\it associative rings with identity, which are algebraic over the field (even Artinian subrings) of own centers} \cite{Pet12}, {\it are stable}.

As in this case for an arbitrary element $a\in R$ and an Artinian subring $K\subseteq \xi R$ the chain of ideals $\left(a\right)\supseteq \left(a^{2} \right)\supseteq \left(a^{3} \right)\supseteq \cdots $ of an Artinian commutative ring $K\left[a\right]$ is stabilized as well. Hence, there exists a positive integer \textit{m} and an element $a'\in R$, which commutes with \textit{a}, such that $a^{m} =a^{m+1} a'$. Then 
\begin{center}
$e=a^{m} \left(a'\right)^{m} =a^{m+1} \left(a'\right)^{m+1} =\cdots =a^{2m} \left(a'\right)^{2m} =e^{2} $  and
\end{center} 
\begin{center}
$a^{m} =a^{m+1} a'=\cdots =a^{2m} \left(a'\right)^{m} =ea^{m}.$
\end{center}

In this case $1-r+ar=a^{m} $. Suppose $r=1+a+\cdots +a^{m-1} $. 

If $g_{1} =g^{-1} g_{0}^{-1} $, $g_2=(g^{-1})^{g_0},$ where $g_{0}=\prod _{l}t_{li} \left(g_{li} r\left(g^{-1} \right)_{ii} \right)$, then
\[\left(g_{1} \right)_{ii} =\left(g_{2} \right)_{ii}=\left(1-r+ar\right)\left(g^{-1} \right)_{ii} =a^{m} \left(g^{-1}
\right)_{ii} \subseteq eR.\]
If $a=\left(g^{-1}\right)_{ii} g_{ii},$ then 
\[\begin{aligned}e&=a^{m} \left(a'\right)^{m}=a^{m+1} \left(a'\right)^{m+1} =a^{m} a\left(a'\right)^{m+1} =a^{m}\left(g^{-1} \right)_{ii} g_{ii} \left(a'\right)^{m+1} \subseteq
\\ &\subseteq\left(g_{1} \right)_{ii} R=\left(g_{2}\right)_{ii} R.
\end{aligned}\]
Thus, Lemma \ref{Petlem11} can be applied to the elements $g_{1} $ and $g_2$. Just as in Theorem \ref{Petthem3} we obtain \textit{g} -- stable element.
\par Let \textit{S} be a multiplicatively closed subset, with identity, of the center $\xi R$ of the ring \textit{R} which does not contain \textit{0, $R_{S}$} -- the classical ring of fractions of the ring \textit{R} by \textit{S}.

The natural homomorphism
\begin{center}
\noindent $\Lambda :R\to R_{S} $, defined by the rule $\Lambda :r\to \frac{r}{1} $
\end{center}
\noindent induces a homomorphism $\Lambda :GL\left(n,R\right)\to GL\left(n,R_{S} \right)$.
\begin{lr}\label{Petlem12} Let \textit{R} be an associative ring with identity, \textit{N} -- a subgroup, invariant under $E(n,R)$ and not containing non-identity transvections. Then $\Lambda N$ does not contain non-identity transvections.
\end{lr}

{\it\textbf{Proof.}} If $\Lambda N$ contains non-identity transvection $\tau$, then there exists a transvection $t\in E_{S} $ such that for some $r \in R,$ $rS\neq0$ the following inclusion holds
\begin{center}
$\Lambda t_{ij} \left(r\right)=\left[\tau ,\Lambda t\right]\in \Lambda N.$
\end{center}
\par This means that $t_{ij}(r)h \in N$ for some $h\in \ker \Lambda $. In such a case there exists $s\in S$, such that $\left(h-1\right)s=0$ and \textit{s} annihilates some non-diagonal element of the matrix $t_{ij}\left(r\right)h$. According to Lemma \ref{Petlem3}, element \textit{s} annihilates all non-diagonal elements of the matrix $t_{ij}\left(r\right)h$. Thus, 
$t_{ij} \left(r\right)s=t_{ij}\left(r\right)hs$ is a diagonal matrix and $rs=0$. The contradiction thus obtained shows that $\Lambda N$  does not contain non-identity transvections.

\begin{lr}\label{Petlem13} Let \textit{R} be an associative ring with identity, \textit{I} -- ideal of \textit{R}, \textit{N} - a subgroup, invariant with respect to $E(n,R)$ and not containing non-identity transvections. If $\Lambda \left(g\right)$ is a stable element of the group $GL\left(n,R_{S} \right)$, then there exists an element $s\in S$ such that for an arbitrary $e\in E\left(n,R\right)$
\begin{center} 
\noindent $\left[g,e,E_{sR} \right]\subseteq E\left(n,I\right)$ if $g\in C_{I} $ and $g\in C\left(n,Anns\right)$ if $g\in N$. \end{center}
\end{lr}

{\it\textbf{Proof.}} Let $g\in C_{I} $. Then $\Lambda \left(g\right)\subseteq C_{I_{S} } $,
\[\left[\Lambda \left(g\right),E\left(n,R_{S} \right)\right]\subseteq E\left(n,I_{S} \right)\]
and there exists an element $s_{0} \in S$ such that $\left[g,e,t_{ij} \left(cs_{0} r\right)\right]\in E\left(n,cI\right)\ker \Lambda $ for arbitrary $e\in E\left(n,R\right)$, $r\in R$, $c\in \xi R$, $1\le i\ne j\le n$.

The inclusion thus obtained holds if the ring \textit{R} should be interchanged with the ring $R\left[x,y\right]$, in which the variables \textit{x} and \textit{y} commute, \textit{x} commutes with the elements of \textit{R}, and \textit{y} - with the elements of $\xi R$. Hence, it can be viewed as a polynomial in terms of variable \textit{x} with the coefficients from the ring $R_{n} \left[y\right]$, which are annihilated by some element $s_{1} \in S$.

Let $s_{ij} =s_{o} s_{1} $. Then $\left[g,e,t_{ij} \left(s_{ij} y\right)\right]\in E\left(n,I\left[y\right]\right)$. Thus, 
\begin{center}
\noindent $\left[g,e,t_{ij} \left(s_{ij} R\right)\right]\subseteq E\left(n,I\right)$     and     $\left[g,e,E_{sR} \right]\subseteq E\left(n,I\right),$
\end{center}
\noindent where $s=\bigcap s_{ij} $ for all pairs $1\le i\ne j\le n$.

Let $g\in N$. In view of Lemma \ref{Petlem12}, the group $\Lambda N$ does not contain non-identity transvections. Since $\Lambda \left(g\right)$ is a stable element, then $\Lambda \left(g\right)\in \xi GL\left(n,R_{S} \right)$. Therefore, there exists $s\in S$ that annihilates non-diagonal element of matrix \textit{g}. In accordance with Lemma \ref{Petlem3}, $g\in C\left(n,Anns\right).$
\begin{sr} Let \textit{R} be an associative ring with identity, $R_{S}$ -- stable rings for all maximal ideals $J_{0}$ of the subring $K\subset \xi R$, $1\in K$, $S=K\backslash J_{0}$. Then \textit{R} is a stable ring.
\end{sr}

{\it\textbf{Proof.}} Let $e_{0} ,e$ be arbitrary elements of the group $E(n, R),$ $G$ -- subgroup of the group $GL(n,R)$ invariant with respect to $E(n,R)$ and $I_{0}$ -- largest ideal of the ring \textit{R} such that $E\left(n,I_{0}\right)\subset G$. Let 
\begin{center}
$J\left(I\right)=\left\{s\in K\left|\left[g,e_{0},e,E_{sR} \right]\subseteq E\left(n,I\right) \,\mbox{if} \,g\in
C\left(n,I\right)\right. \right\}$ and 
\end{center}
\begin{center}
$J\left(G\right)=\left\{s\in K\left|{\rm \; }\Lambda _{0}\left(g\right)\in C\left(n,Ann\Lambda _{I_{0} }
\left(s\right)\right)\, \mbox{if} \,g\in G\right. \right\}.$
\end{center}
It is understandable that $J\left(I\right)$ and $J\left(G\right)$ -- ideals of the ring \textit{K}.

If $J\left(I\right)\ne K$, then there exists a maximal ideal $J_{0}\left(I\right)$ of the ring \textit{K} such that 
\begin{center}
$J\left(I\right)\subseteq J_{0} \left(I\right), S=K\backslash J_{0}\left(I\right).$
\end{center} 
Similarly one can define $S=K\backslash J_{0} \left(G\right)$, if $J\left(G\right)\ne K$, where $J\left(G\right)\subseteq J_{0}
\left(G\right)$.

Let $g\in C\left(n,I\right)$. Since $R_{S}$ is a commutator ring, then $\Lambda _{I_{0} } \left[g,e_{0} \right]\subseteq
E\left(n,I_{S} \right)$. According to Lemma \ref{Petlem13}, there exists an element in \textit{S}, which is contained in $J\left(I\right)$. The contradiction thus obtained shows that $J\left(I\right)=K$, $1\in J\left(I\right)$, \textit{R} -- weakly-commutator ring.

Let $g \in G$, $\overline{R}={\raise0.7ex\hbox{$ R$}\!\mathord{\left/{\vphantom{R I_{0}}}\right.\kern-\nulldelimiterspace}\!\lower0.7ex\hbox{$ I_{0}  $}}$. As in Lemma \ref{Petlem12} we prove that $\Lambda_{I_{0}}(G)$ does not contain non-identity transvections. Since $I_{0} \bigcap K\subseteq J\left(G\right)$, then $\overline{S}$ does not contain zero element. As a quotient of the normal ring $R_{S}$, the ring $\overline{R}_{\overline{S}} $ is partially normal.

If $\Lambda :\overline{R}\to \overline{R}_{\overline{S}} $, then the group $\Lambda \Lambda _{I_{0} } \left(G\right)$, according to Lemma \ref{Petlem12}, does not contain non-identity transvections. In view of the partial normality of the ring $\overline{R}_{\overline{S}} $ and Lemma \ref{Petlem13} there exists an element in \textit{S}, which is contained in $J\left(G\right)$. Therefore, $J\left(G\right)=K$, $1\in J\left(G\right)$, $\Lambda _{I_{0} } \left(g\right)\in \xi
GL\left(n,R\right)$, $g\in C\left(n,I_{0} \right)$. This means that \textit{R} is a normal and, as a consequence, stable ring.

\par In the particular case, when $R_{S} $ - {\it rings which satisfy the condition of stability of rank $n-1>1$ for all maximal ideals {\it J} of the subring $K\subset \xi R$ under the condition $1\in K$, $S=K\backslash J$ the stability of the ring $R$ is proved in} \cite{Pet9}.

\begin{sr}\label{Petnas7} \cite{Pet27} Let {\it R} be an associative ring with identity, which is integer-algebraic over the subring $K\subset \xi R$, $1\in K$. Then {\it R} is a stable ring.
\end{sr}

{\it\textbf{Proof.}} Let $I_{0}$ be a maximal ideal of the ring {\it K}, $S=K\backslash I_{0}, r\in R.$ Then $K_{S} \left(r\right)$ is a finitely-generated module over $K_{S} $. According to Nakayama's Lemma $J\left(K_{S} \right)\subseteq J\left(K_{S} \left(r\right)\right)$ and, as a consequence, $J\left(K_{S} \right)\subseteq J\left(R_{S} \right)$. Therefore, 
${\raise0.7ex\hbox{$ R_{S}  $}\!\mathord{\left/{\vphantom{R_{S}  J\left(R_{S} \right)}}\right.\kern-\nulldelimiterspace}\!\lower0.7ex\hbox{$ J\left(R_{S} \right) $}} $~- is a ring, algebraic over the field ${\raise0.7ex\hbox{$ K_{S}  $}\!\mathord{\left/{\vphantom{K_{S}  J\left(K_{S} \right)}}\right.\kern-\nulldelimiterspace}\!\lower0.7ex\hbox{$ J\left(K_{S} \right) $}} $.
In such a case, as was mentioned above, the ring ${\raise0.7ex\hbox{$ R_{S}  $}\!\mathord{\left/{\vphantom{R_{S}  J\left(R_{S} \right)}}\right.\kern-\nulldelimiterspace}\!\lower0.7ex\hbox{$ J\left(R_{S} \right) $}} $ is stable. In view of Corollary \ref{Petnas5}, $R_{S} $ is a stable ring and, according to the corollary 6, {\it R} is a stable ring.

As is known \cite{Pet16} not every associative ring with identity is stable. For instance, {\it algebra over field with $2n^{2} $ generating elements $x_{ij} ,y_{ij} $,  $1\le i,j\le n$ and the defining relations $$\left(x_{ij} \right)\left(y_{ij} \right)=\left(y_{ij} \right)\left(x_{ij} \right)=1$$ is not a stable ring}.

\par However, the class of stable rings is quite wide. The most vividly it was demonstrated in the work \cite{Pet18}.
\begin{vr}\label{Petozn15} Let {\it R} be an associative ring with identity. Ideal {\it F} of the ring {\it R} is called weakly Noetherian (respectively integer-algebraic) if for arbitrary elements $y,z\in F$ ,  $m\ge 1$ left and right modules 
\begin{center}
 $\sum \limits _{m}Rzy^{m}$  and  $\sum_m y^mzR$      (respectively  $\sum \limits _{m} \xi Rzy^m $ and $\sum \limits _{m} \xi Ry^mz $ )
\end{center}
are finitely generated as modules over {\it R} (over $\xi R$ respectively).
\end{vr}
\begin{vr} Associative ring {\it R} with identity is called weakly Noetherian (respectively integer-algebraic) if there exists a chain of ideals $$0=I_{0} \subseteq I_{1} \subseteq \ldots \subseteq I_{q+1} =R$$ such that the ideals ${\raise0.7ex\hbox{$ I_{i+1} $}\!\mathord{\left/{\vphantom{I_{i+1}  I_{i}}}\right.\kern-\nulldelimiterspace}\!\lower0.7ex\hbox{$ I_{i}  $}} $ in the rings  ${\raise0.7ex\hbox{$ R$}\!\mathord{\left/{\vphantom{R I_{i}}}\right.\kern-\nulldelimiterspace}\!\lower0.7ex\hbox{$ I_{i}  $}}$ are weakly Noetherian (integer-algebraic respectively) for all $1\le i\le q$.
\end{vr}

Obviously, the block integer-algebraic rings are weakly Noetherian.
\par It is known \cite{Pet21, Pet22, Pet10} that {\it PI} -- rings that are block integer-algebraic. Obviously, rings that are algebraic over subrings of own centers are block integer-algebraic.

Let $g\in GL\left(n,R\right)$  and {\it l} - maximal integer such that $I_{l} g_{1n} =0$. If $l<q+1$, then we choose $g_{1} \in \left[g,t_{n1} \left(I_{l+1} \right)\right]$,  $y=\left(g_{1} \right)_{11}$,  $z=\left(g_{1} \right)_{1n}$. Then $y-1$   and {\it z} are contained in $I_{l+1} \bigcap g_{1n} R$ and $I_{l} \left(y-1\right)=0$.

Therefore, there exists a positive integer {\it m} such that 
\begin{center}
 $z\left(y-1\right)^{m} -\sum \limits _{p=1}^{m-1}s_{p} z\left(y-1\right)^{p} \in I_{l}$, $zy^{m+1} =\sum \limits _{p=0}^{m}r_{p} zy^{p}$, $r_{p} ,s_{p} \in R$.
\end{center}

\par Let {\it N} be a subgroup of $GL\left(n,R\right)$, invariant with respect to $E\left(n,R\right)$ and not containing non-identity transvections. If $g\in N$, then  $g_{1} \in N$ and, according to Lemma \ref{Petlem3}, {\it y} is not a divisor of zero. According to Lemma \ref{Petlem4}, from the equality $ry+r_{0} z=0$, where {\it r} - some element of the ring {\it R}, it follows that $r_{0}z=r=0.$ Similarly one can prove $0=r_{0} z=\cdots =r_{m} z=z.$ Hence, $\left[g,t_{n1} \left(I_{l+1} \right)\right]\subseteq \xi GL\left(n,R\right).$ In view of Lemma \ref{Petlem3}, $I_{l+1}g_{1n} =0$. The contradiction thus obtained shows that $l=q+1$, $g_{1n} =0$, {\it R} - partially normal ring.

Therefore, {\it weakly Noetherian rings are partially normal}. Since the property of being weakly Noetherian is preserved under factorization, then {\it all the quotients of weakly Noetherian rings are partially normal as well}.

Let {\it I} be an ideal of {\it R}, which is contained in some weakly Noetherian ideal of the ring {\it R}, {\it y} - arbitrary element of {\it I}. Then $zy\in I$ and there exists a positive integer {\it m} such that 
\begin{center}
 $zy^{m+1} =\sum \limits _{p}r_{p} zy^{p}$, where  $z,r_{p} \in R$ ,  $1\le p\le m.$
\end{center}
\par Let $\lambda \in \xi R.$ Multiply the equality above by $\lambda ^{m+1} $. Since $\lambda y =\lambda y -1+1,$ then there exists a polynomial $\psi \left(\lambda \right)$ such that the following equality holds 
\begin{center}
 $\psi \left(\lambda \right)z+a\left(1-\lambda y\right)=0$ , where  $\psi \left(0\right)=1$ .
\end{center}
\par Let $g\in C\left(n,I\right)$,  $g_{\lambda } =t_{pq} \left(\lambda r\right)^{g}$, $c_{\lambda } =\left[t_{pq} \left(-\lambda r\right),g\right]$,  $y=1-\left(g_{1} \right)_{ii}$, \\$z=\left(g_{1}^{-1} \right)_{jk} \left(g_{1} \right)_{ki},$ where $i,j,k$ - pairwise distinct numbers, $r\in R$.

Then $g_{\lambda } =\lambda g_{1} -\lambda +1,$  $c_{\lambda } \in C_{I},$   $g_{\lambda } =t_{pq} \left(\lambda r\right)c_{\lambda }.$

\par Since 
\[
(g_\lambda)_{ii}=\lambda(g_1)_{ii}-\lambda+1=1-\lambda y,\; (g_\lambda)_{ki}=\lambda(g_1)_{ki},\]\[ (g_\lambda^{-1})_{jk}=(g_{-\lambda})_{jk}=-\lambda (g_1)_{jk}=\lambda (g_1^{-1})_{jk} \; \mbox{and}\; t_{pq} \left(-\lambda r\right)_{jk}t_{pq} \left(\lambda r\right)_{ki}=0,
\]
then
\begin{center}
     $\psi \left(\lambda \right)\left(g_{\lambda }^{-1} \right)_{jk} \left(g_{\lambda } \right)_{ki} +\lambda ^{2} a\left(g_{\lambda } \right)_{ii}      =0,$ $a \in I.$
     \end{center}
In accordance with Lemma \ref{Petlem7},
\begin{center}
$\left[c_{\lambda } ,t_{ij} \left(J\psi _{ij} \left(\lambda \right)\right)\right]\subset E\left(n,I\right)\bigcap
E\left(n,J\right),$
\end{center}
where $J$ is an arbitrary deal of $R.$

\par Following similar ``right-sided'' arguments and taking into account the matrix commutator formulas for each pair there exists a polynomial such that 
\begin{center}
 $\left[c_{\lambda } ,t_{ij} \left(J\psi _{ij} \left(\lambda \right)J\right)\right]\subset E\left(n,I\right)\bigcap E\left(n,J\right)$, where  $\psi _{ij} \left(0\right)=1$.
\end{center}
Let $f\left(\lambda \right)=\prod \limits _{i,j}\psi _{ij} \left(\lambda \right) $ for all pairs $1 \leq i \neq j \leq n.$ Then 
\[\left[c_{\lambda } ,E_{Jf\left(\lambda \right)J}\right]\subseteq E\left(n,I\right)\bigcap E\left(n,J\right).\]

\par Define $I_{1} =Rf\left(\lambda \right)R$, $I_{2} =I_{1}^{2} f\left(1-\lambda \right)I_{1}^{2}$.

Obviously $I_{1}, I_{2} $ are ideals of the ring {\it R} and 
\begin{center} 
$\left[c_{\lambda } ,E_{I_{1} } \right]\subseteq E\left(n,I\right),$       $\left[c_{1-\lambda } ,E_{I_{2} } \right]\subseteq E\left(n,I\right).$ 
\end{center}
Since $E_{I_{2} } \subseteq E\left(n,I_{2} \right)\subseteq E\left(n,I_{1}^{2} \right)\subseteq E_{I_{1} } $, then for an arbitrary element \\$e\in E\left(n,R\right)$ the following inclusions hold 
\begin{center}
 $\left[c_{\lambda }^{e} ,E_{I_{1} } \right]\subseteq \left[c_{\lambda } ,E\left(n,I_{2} \right)\right]^{e} \subseteq E\left(n,I\right)$     and     $\left[c_{\lambda } ,E_{I_{1} } \right]^{e} \subseteq E\left(n,I_{I}^{2} \right)\subseteq E_{I_{1} }$.
\end{center}
Taking into account that $c_{1} =c_{\lambda }^{t_{pq} \left(\lambda -1\right)r} c_{1-\lambda } $ we obtain $\left[c_{1} ,E_{I_{2} } \right]\subseteq E\left(n,I\right)$.

Let us put $I_{0} =\sum \limits _{\lambda }I_{2}$ for all $\lambda \in \xi R$. If $I_{0} \ne R$ then, due to the equality $f\left(0\right)=1$, the image of the polynomial 
\begin{center}
$f\left(\lambda \right)^{2} f\left(1-\lambda \right)f\left(\lambda \right)^{2} $
\end{center}
 is non-zero over the ring ${\raise0.7ex\hbox{$ R $}\!\mathord{\left/{\vphantom{R I_{0}
}}\right.\kern-\nulldelimiterspace}\!\lower0.7ex\hbox{$ I_{0}  $}}$, and images of the elements of $\xi R$ are its roots.

If $\xi R$ contains an infinite field, then $I_{0} =R$,  $\left[c_{1} ,E\left(n,R\right)\right]\subseteq E\left(n,I\right)$, {\it R} - weakly-commutator ring.

In accordance with Lemma \ref{Petlem2},  {\it R} is a stable ring. Thus, the {\it weakly Noetherian rings, which contain infinite fields in own centers, are stable} \cite{Pet18}.

In the particular case, the {\it block integer-algebraic rings are stable} without the demand of existence of infinite fields in own centers.

Indeed, if elements of the ideal {\it I} of the ring {\it R} are integer-algebraic over the subring $K\subseteq \xi R$,  $1\in K$,  $r\in I$,  $I_{0}$ is a maximal ideal of {\it K},  $S=K\backslash I_{0}$, then  $K_{S} \left(r\right)$ is a finitely generated module over $K_{S} $. Due to Nakayama's Lemma,
\begin{center}
 $J\left(K_{S} \right)\subseteq J\left(K_{S} \left(r\right)\right)$         and          $J\left(K_{S} \right)\subseteq J\left(K_{S} \left(I\right)\right)$.
\end{center}
Since the ring ${\raise0.7ex\hbox{$ K_{S} \left(I\right) $}\!\mathord{\left/{\vphantom{K_{S} \left(I\right) J\left(K_{S} \left(I\right)\right)}}\right.\kern-\nulldelimiterspace}\!\lower0.7ex\hbox{$ J\left(K_{S} \left(I\right)\right) $}} $ is algebraic over the field ${\raise0.7ex\hbox{$ K_{S}  $}\!\mathord{\left/{\vphantom{K_{S}  J\left(K_{S}\right)}}\right.\kern-\nulldelimiterspace}\!\lower0.7ex\hbox{$ J\left(K_{S} \right) $}} \cong {\raise0.7ex\hbox{$ K $}\!\mathord{\left/{\vphantom{K I_{0} }}\right.\kern-\nulldelimiterspace}\!\lower0.7ex\hbox{$ I_{0}  $}},$ the rings $K_{S} \left(I\right)$ and, respectively, $K\left(I\right)$ are stable. In such a case the groups $\left[C_{I} ,E_{K} \right]$  and  $\left[C_{I} ,E\left(n,R\right)\right]$ are contained in $E\left(n,I\right)$.

Hence, if {\it R} is a block integer-algebraic ring then, according to the fact proved above,

\begin{center}
$\left[C_{I_{i+1} } ,E\left(n,R\right)\right]\subset E\left(n,I_{i+1} \right)C_{I_{i} } $
\end{center}
for all $0\le i\le q.$ Therefore, {\it R} is a weakly-commutator ring with partially normal quotients. In view of Lemma \ref{Petlem2}, $R$ is a stable ring \cite{Pet19,Pet28}.


{

\end{document}
\begin{thebibliography}{99}
\bibitem{Pet20}
E. Abe, {\it Automorphisms of Chevalley groups over commutative rings,} Algebra and calculus. {\bf 5} (1993), no. 2, pp. 74--90.
\bibitem{Pet22}
E. Abe, {\it Chevalley groups over commutative rings,} Proc.-Conf.-Radical-Theory.-Sendai. (1998), pp. 1--23.
\bibitem{Pet21}
E. Abe, K. Suzuki, {\it On normal subgroups of Chevalley groups over commutative rings,} T\^{o}hoku Math.J. {\bf 28} (1976), no. 1, pp. 185--198.
\bibitem{Pet26}
V.A. Artamonov, {\it Serre's quantum problem,} Successes of mathematical sciences. {\bf 53} (1998), no. 4(322), pp. 3--76.
\bibitem{Pet1}
H. Bass, {\it K-theory and stable algebra,} Publ.Math.IHES. {\bf 22} (1964), pp. 5--60.
\bibitem{Pet2}
H. Bass, {\it Algebraic K-theory,} Mir, Moscow, (1973), 591 p.
\bibitem{Pet11} 
I.Z. Borevich, N.A. Vavilov, {\it Location of subgroups in the full linear group over the comutative ring,} Tr. MIAN USSR {\bf 165} (1984), pp. 24--42.
\bibitem{Pet16}
V.I. Gerasimov, {\it Group of units of the free product of rings,} Math.col. {\bf 134} (1987), no. 1, pp. 42--45.
\bibitem{Golubchik}
I. Z. Golubchik, {\it On a full linear group over associative rings,} UMN {\bf 27} (1973), no. 3, pp. 179--180.
\bibitem{Pet10}
I.Z. Golubchik, {\it On subgroups of a general linear group $GL\left(n,R\right)$ over associative ring {\it R},} UMN. {\bf 39} (1984), no. 1, pp. 125--126.
\bibitem{Pet8}
I.Z. Golubchik, {\it Isomorphism of the General Linear Group $GL_n(R), n\geq 4$ over an associative Ring,} Contemporary
Mathematics. {\bf 131} (1992), part 1, pp. 123--136.
\bibitem{Pet29}
I.Z. Golubchik, {\it Isomorphisms of projective groups over associative rings,} Fundament. and appl. math. (1995), pp. 311–-314.
\bibitem{Pet18}
I.Z. Golubchik, {\it On a full linear group over weakly Noetherian associative algebras,} Fundament. and applied mathematics {\bf 1} (1995), no. 3, pp. 661--668.
\bibitem{Pet19}
I.Z. Golubchik, {\it On a full linear group over block-algebraic rings,} Thesis dokl.inter.algebr.conf.in memory of Fadeev D.K. San-Petersburg (1997), pp. 186--187.
\bibitem{Pet25}
I.Z. Golubchik, {\it Lie-type groups over PI-rings,} Fundamental and applied mathematics. {\bf 3} (1997), no. 2, pp. 399--424.
\bibitem{Pet30}
I.Z. Golubchik, A.V. Mikhalev, {\it Isomorphism of a complete linear group over an associative ring,} Vestnik Moscov. Univ.
Series 1. Math. Mech. (1983), no. 3, pp. 61--72.
\bibitem{Pet15}
A.I. Hahn and O.T. O'Meara, {\it The Classical Groups and K-theory,} Springer -Verlag. (1989), 494 p.
\bibitem{Pet13}
S.H. Khlebutin, {\it Sufficient conditions for normality of the subgroup of elementary matrices,} UMN. {\bf 39} (1984), no. 3, pp. 245--246.
\bibitem{Pet12}
S.H. Khlebutin, {\it Some properties of elementary subgroup. Algebra: logic and number theory,} Izd MSU, Moscow (1986), pp. 86--90.
\bibitem{Pet14}
V.M. Petechuk, {\it Homomorphisms of linear groups over commutative rings,} Mathematical notes. {\bf 46} (1989), no. 5, pp. 50--61.
\bibitem{Pet28}
V.M. Petechuk, {\it Stable structure of linear groups over rings,} Dopovidi NAN of Ukraine (2001), no. 11, pp. 17--22.
\bibitem{Pet27}
V.M. Petechuk, {\it Stability structure of linear group over rings,} Matematychni studii. {\bf 16} (2001), no. 1, pp. 13--24.
\bibitem{Uzhg}
V.M. Petechuk, {\it Stability of rings,} Scientific notices of Uzhgorod University, mathematics and informatics section (2009), no. 19, pp. 87--111.
\bibitem{Pet3}
D.A. Suprunenko, {\it Groups of matrices,} Nauka, Moscow, (1972), 351 p.
\bibitem{Pet6}
A.A. Suslin {\it On the structure of special linear group over the ring of polynomials,} Izv. AN USSR. Ser.math. {\bf 41} (1977), no. 2, pp. 235--252.
\bibitem{Pet23}
G. Taddei, {\it Normalite des groupes elementaire class les groupes de Chevalley sur anneau,} Cont.-Math.-Amer.-Math.-Soc. {\bf 55} (1986), part II, pp. 693--710.
\bibitem{Pet4}
L.N. Vaserstein, {\it On stabilization in algebraic K-theory,} Func. analysis and its applications. {\bf 3} (1969), no. 2, pp. 85--86.
\bibitem{Pet5}
L.N. Vaserstein, A.A. Suslin, {\it Serre's problem on projective modules over polynomial rings and algebraic K-theory,} Math.USSR Izv. {\bf 10} (1976), no. 5, pp. 937--1001.
\bibitem{Pet9}
L.N. Vaserstein, {\it On the normal subgroups  $GL_{n} $  over a ring,} Springer Lecture Notes. 854 (Algebraic K-theory) (1981), pp. 454--465.
\bibitem{Pet17}
L.N. Vaserstein, {\it Normal subgroups of general linear groups over von Neumann regular rings,} Proc.Amer.Math.Soc. {\bf 96} (1986), no. 2, pp. 209--214.
\bibitem{Pet24}
L.N. Vaserstein, {\it On normal subgroups of Chevalley groups over commutative rings,} T\^{o}hoku Math.J. {\bf 38} (1986), pp. 219--230.
\bibitem{Pet7}
J.S. Wilson, {\it The normal and subnormal structure of general linear groups,} Proc. Cambr. Phil. Soc. {\bf 71} (1972), no. 2, pp. 163--177.

\end{thebibliography}
